\documentclass[12pt]{article}
\usepackage{latexsym,amsfonts,amsmath,graphics,amssymb}
\usepackage{epsfig}
\usepackage{verbatim}
\usepackage{stix}

\usepackage{pgfplots}%Fuer farbige Kommentare...
\usepackage{mathtools}
\usepackage{esint}

\DeclareMathOperator{\supp}{supp}
\DeclareMathOperator{\wor}{wor}

\DeclareMathOperator{\dom}{dom}

\newtheorem{theorem}{Theorem}
\newtheorem{lemma}{Lemma}
\newtheorem{corollary}{Corollary}
\newtheorem{proposition}{Proposition}

\newtheorem{definition}{Definition}

\newtheorem{remark}{Remark}
\newtheorem{example}{Example}
\newenvironment{proof}{\begin{trivlist}
    \item[\hskip\labelsep{\it Proof.}]}{$\hfill\Box$\end{trivlist}}

\newcommand{\bsxi}{\boldsymbol{\xi}}

\newcommand{\bsx}{\boldsymbol{x}}

\newcommand{\bst}{\boldsymbol{t}}

\newcommand{\bsp}{\boldsymbol{p}}
\newcommand{\bsone}{\boldsymbol{1}}

\newcommand{\bsy}{\boldsymbol{y}}

\newcommand{\bsalpha}{\boldsymbol{\alpha}}

\newcommand{\bszero}{\boldsymbol{0}}
\newcommand{\rd}{\,\mathrm{d}}
\newcommand{\NN}{\mathbb{N}}

\newcommand{\RR}{\mathbb{R}}
\newcommand{\real}{\mathbb{R}}

\newcommand{\R}{\mathbb{R}}         % Dieser Block kommt aus MA_R+  

\newcommand{\F}{\mathcal{F}}
\newcommand{\Sc}{\mathcal{S}}
\DeclareMathOperator{\exto}{ext_0}
\DeclareMathOperator{\exti}{ext_1}

\DeclareMathOperator{\sgn}{sgn}
\DeclareMathOperator{\osc}{osc}
\DeclareMathOperator{\Osc}{osc_+^1}

\newcommand{\HH}{\mathcal{H}}

\newcommand{\cf}{{\mathcal{F}}}
\newcommand{\cfi}{{\mathcal{F}}^{-1}}

\newcommand{\sign}{{\rm sign} \, }
\newcommand{\N}{\mathbb{N}}

\newcommand{\be}{\begin{equation}}
\newcommand{\ee}{\end{equation}}
\newcommand{\beq}{\begin{eqnarray}}
\newcommand{\beqq}{\begin{eqnarray*}}
\newcommand{\eeq}{\end{eqnarray}}
\newcommand{\eeqq}{\end{eqnarray*}}

\newcommand{\C}{\mathbb{C}}
\allowdisplaybreaks

\newcommand{\jo}[1]{\textcolor{purple}{#1}}

\title{On Sobolev spaces and Koksma-Hlawka-type inequalities of fractional order}

\author{
$\quad$Yannick Meiners\\$\quad$\\
{\small Institut f\"ur Mathematik,
Universit\"at Osnabr\"uck}\\
{\small Albrechtstra{\ss}e 28 a} \\
{\small 49076 Osnabr\"uck, Germany} \\
{\small email: yannick.meiners@uni-osnabrueck.de}\\
$\quad$\\
Michael Gnewuch\\$\quad$\\
{\small Institut f\"ur Mathematik,
Universit\"at Osnabr\"uck}\\
{\small Albrechtstra{\ss}e 28 a} \\
{\small 49076 Osnabr\"uck, Germany} \\
{\small email: michael.gnewuch@uni-osnabrueck.de}\\
}

\date{\today}

\begin{document}

\maketitle

\pagebreak

\begin{comment}
    Key words: Bessel potential spaces, fractional Sobolev spaces, fractional derivative, Riemann-Liouville fractional derivative (?), Koksma-Hlawka inequality, QMC-Quadrature,
    
\end{comment}

\begin{abstract}
\noindent
Spaces of fractional smoothness in the sense of Riemann-Liouville of order $0<\alpha<1$ were introduced in [J. Dick, Ann. Math. Pura Appl., 187, 385-403, 2008]. 
 Dick also proved Koksma-Hlawka type inequalities for these spaces, using a notion of fractional discrepancy.\\
In [M. Gnewuch et al., Constr. Approx., 63, 647-711, 2026]
upper error-bounds for QMC-integration on these spaces were obtained using $(t,m,s)$-nets and it was proved that the above mentioned Koksma-Hlawka inequalities by Dick are sharp. One step in the proof of the upper bounds was the embedding of the spaces of fractional smoothness into 
Bessel potential spaces, also called fractional Sobolev spaces.\\
We prove that, up to norm equivalence, spaces of fractional smoothness and Bessel potential spaces are equal. From this we deduce that the upper error bounds for QMC-integration from Gnewuch et al. are indeed optimal for spaces of fractional smoothness. Furthermore, we obtain two types of Koksma-Hlawka inequalities for Bessel potential spaces. 
The first type uses the fractional discrepancy as introduced by Dick.
The second type adopts an approach by Triebel from [H. Triebel, EMS Press, 2010] by using the norm of the (local) discrepancy function in an adjoint Bessel potential space.
%
%the second type uses the norm of a fractional Sobolev space on the discrepancy function, similar to the Koksma-Hlawka type inequalities proved in [H. Triebel, EMS Press, 2010] for Besov spaces.
%We also obtain Koksma-Hlawka type inequalities for Bessel potential spaces in which the discrepancy is also measured in the norm of a Bessel potential space, similar to the Koksma-Hlawka type inequalities proved in [Triebel, EMS Press, 2010] for Besov spaces.
\end{abstract}


%&&&&&&&&&&&&&&&&&&&&&&&&&&&&&&&&&&&&&&&&&&&&&&&&&&&&&&

\tableofcontents

\section{Introduction}

\subsection{Goal and Motivation}\label{GandM}

Consider the integration problem of approximating an $s$-dimensional integral of the form
\begin{equation}\label{eq:I_s}
I_s(f) := \int_{[0,1]^s} f(\bsx) \, {\rm}d \bsx
\end{equation}
with a \emph{quasi-Monte Carlo (QMC) quadrature rule}, i.e., an equal weight quadrature rule whose weights sum up to one. The classical Koksma-Hlawka inequality 
ensures that the integration error of a QMC quadrature rule 
\begin{equation}\label{eq:qmc_cub}
Q_P(f) := \frac{1}{|P|} \sum_{\bsp\in P} f(\bsp),
\end{equation}
where $|P|$ denotes the cardinality of the set of integration points $P$, 
can be bounded by the product of the star discrepancy of $P$ and the Hardy-Krause variation of the integrand $f$, see \cite{ Kok42,Hla61}. Unfortunately, many integrands that appear in important applications do not exhibit a finite Hardy-Krause variation, see, e.g., \cite{Owe05}.

To overcome this fundamental issue, %our co-author 
Dick~\cite{Dic08a} studied the  integration problem for a suitable range of parameters $p,q \in [1,\infty]$ on Banach spaces $\mathcal{H}_{\alpha, s, p,q}$ of functions $f: [0,1]^s\to \RR$ with fractional smoothness $\alpha \in (0,1)$, defined via the fractional derivative in the sense of Riemann and Liouville. These spaces contain functions that do not have finite Hardy-Krause variation. 
%In \cite{D08} 
He introduced a notion of fractional discrepancy $D^*_{\alpha,s,p,q}$ and derived a Koksma-Hlawka-type inequality of fractional order. It ensures that the integration error 
of $Q_P$ 
can be bounded by the product of a fractional discrepancy of $P$ and the norm $||\cdot||_{\alpha,s,p,q}$ of the integrand in
$\mathcal{H}_{\alpha, s, p,q}$,
\begin{equation}\label{ineq:FracKoksma} 
|I_s(f) - Q_P(f)| \le  D^*_{\alpha,s,p',q'}(P) \|f\|_{\alpha,s,p,q}
\hspace{3ex}\text{for all $f\in \mathcal{H}_{\alpha, s, p,q}$\,,}
\end{equation}
where $p'$ and $q'$ denote the adjoint exponents of $p$ and $q$, respectively, i.e.,
$$
\frac{1}{p} + \frac{1}{p'} = 1 = \frac{1}{q} + \frac{1}{q'}\,.
$$
In \cite[Thm.~6]{GDMS24} it was verified that this inequality is actually sharp. More precisely, if one 
considers the worst-case error,  
i.e., the supremum of the error on the left-hand side in \eqref{ineq:FracKoksma} taken over all functions $f$ from the norm unit ball in $\mathcal{H}_{\alpha, s, p,q}$, then we  obtain the identity
\begin{equation}
\label{id:FracKoksma}
e^{\wor}(Q_P,\HH_{\alpha,s,p,q}) =  D^\ast_{\alpha,s,p',q'}(P).
\end{equation}
Also the question of what kind of integration nodes actually have a small fractional discrepancy was (partially) answered in that paper---it turns out that 
(arbitrary) $(t,m,s)$-nets $P$ in base $b$ satisfy the bound
\begin{equation}\label{frac_disc_estimate2}
%e^{\wor}(Q_P, \HH_{\alpha,s,p, q}) = 
D^*_{\alpha,s,p',q'}(P) 
\le C N^{-\alpha} \ln(N)^{\frac{s-1}{2}}, 
\end{equation}
see \cite[Cor.~8]{GDMS24}. 
The following question however remained unanswered.\\

\noindent
\begin{itemize}
     \item[(Q1)]\textit{Is the bound \eqref{frac_disc_estimate2} on the fractional discrepancy best possible or 
     are there Quasi-Monte Carlo
point sets that exhibit a substantially smaller discrepancy?}
\end{itemize}

%\noindent
In the course of the proof of \eqref{frac_disc_estimate2} it was shown that the fractional spaces $\mathcal{H}_{\alpha, s, p,q}$ 
are continuously embedded in Bessel potential spaces $S^\alpha_pH([0,1]^s)$, which are also referred to as fractional Sobolev spaces in the literature,
see \cite[Cor.~3]{GDMS24}. This gave rise to another question.\\

\noindent
\begin{itemize}
     \item[(Q2)]\textit{Are fractional spaces in the sense of Riemann-Liouville and Bessel potential spaces even equal up to norm-equivalence?}\\
\end{itemize}

%In this paper we want to investigate these two related questions and present the complete answer.
Furthermore we are interested in bounds for the integration error of functions from Bessel potential spaces. The results from \cite{Dic08a, GDMS24} motivate the final question of this paper.\\
\begin{comment}
Furthermore we are interested in bounds for the integration error.\\
For a pointset $P=\{\bsx_0,\ldots,\bsx_{N-1}\}$ in $[0,1]^s$ with $\bsx_n = (x_{n,1},\ldots,x_{n,s})$
and $\emptyset\neq u\subseteq S$ we define the \emph{(local) discrepancy function $\Delta$} for $\bst_u\in [0,1]^{u}$ as
\begin{equation}
   \Delta(\bst_u,u,P):=\prod_{j\in u}t_j-\frac 1N\sum_{n=0}^{N-1}\prod_{j\in u}\mathbb 1_{( x_{n,j},1]}(t_j); 
\end{equation}
in case $u = S$, we usually suppress the references to $u$ and obtain 
the classical, well-studied discrepancy function $\Delta(\cdot, P)$, anchored in $\bszero$. The \emph{($L^\infty$-)star discrepancy
of $P$} is defined as 
\[
D^*_{N}(P) := \sup_{\bsx\in [0,1)^s}|\Delta(\bsx,P)|
\]
%Let us use the short hands
%\begin{equation*}
%I_s(f):= \int_{[0,1]^s} f(\bst)\,{\rm d}\bst
%\hspace{3ex}\text{and}\hspace{3ex}
%Q_P(f) := \frac{1}{N}
%\sum_{n=0}^{N-1} f(\bsx_n).
%\end{equation*}
and the \emph{classical Koksma-Hlawka inequality}, see
\cite{ Kok42, Hla61}, is of the form
\begin{equation}\label{KHII}
\left|I_s(f) - Q_P(f) \right|  \le  D^\ast_{N}(P)
V_{{\rm HK}}(f),
\end{equation}
where $V_{{\rm HK}}(f)$ is the \emph{variation of $f$ in the sense of Hardy and Krause}.
This motivates the final question of this paper.\\
    
\end{comment}

\noindent
\begin{itemize}
     \item[(Q3)]\textit{Are there Koksma-Hlawka-type inequalities similar to \eqref{ineq:FracKoksma} for functions $f\in S^\alpha_pH([0,1]^s)$?}\\
\end{itemize}
In particular, we are interested in Koksma-Hlawka-type inequalities for Bessel potential spaces in terms of the fractional discrepancy as introduced by Dick, as well as in terms of the norm of the local discrepancy function taken in a suitable Bessel potential space, analogously to the Koksma-Hlawka-type inequalities deduced in \cite{Tri10} for Besov spaces.\\

\noindent 
In this paper we want to investigate these three related questions and present complete answers.

\begin{comment}
A question that remains is whether this bound on the fractional discrepancy is best possible or if we can find Quasi-Monte Carlo
point sets that exhibit a substantially smaller discrepany?

In the course of the proof of \eqref{frac_disc_estimate2} it was shown that the fractional spaces $\mathcal{H}_{\alpha, s, p,q}$ 
are continuously embedded in Bessel potential spaces $S^\alpha_pH([0,1)^s)$, which are also referred to as fractional Sobolev spaces in the literature,
see \cite[Cor.~3]{GDMS24}.
But the question whether both spaces are maybe even equal (up to norm equivalence) was left open in \cite{GDMS24}.\\

In this paper we want to investigate these two related questions and present the complete answer, as well as introduce a Koksma-Hlawka inequality for Bessel potential spaces (also known as fractional Sobolev spaces).
\end{comment}

\subsection{Main Results}

To be more precise, we will establish the following results:\\
After discussing the univariate spaces of fractional smoothness in the sense of Riemann-Liouville $\mathcal H_{\alpha,p}:=\mathcal H_{\alpha,1,p,q}$ (these spaces are norm-equivalent for different values of $q$) and the univariate fractional Sobolev spaces $H^\alpha_p(\R)$ and $H^\alpha_p([0,1])$, we conclude Section 2 by proving the following embedding result. Let $1<\alpha^{-1}<p<\infty$. Then
$$H^\alpha_p([0,1])\hookrightarrow\mathcal H_{\alpha,p},$$
see Theorem \ref{embedding01}. Since the inverse embedding has been proved earlier in \cite{GDMS24}, we get as a corollary 
$$H^\alpha_p([0,1])=\mathcal H_{\alpha,p}$$
up to equivalence of norms, see Corollary \ref{equal01}.

In Section 3 we transfer these results to the multivariate case. For $1<\alpha^{-1}<p<\infty$ we get, up to equivalence of norms, the identity
\begin{equation}
\label{idSpaces}
S^\alpha_p H([0,1]^s)= \mathcal H_{\alpha ,s,p },
\end{equation}
where $S^\alpha_p H([0,1]^s)$ are Bessel potential spaces or fractional Sobolev spaces of dominated mixed smoothness and $\mathcal H_{\alpha ,s,p }$ are multivariate spaces of fractional smoothness in the sense of Riemann-Liouville, see Corollary \ref{equal02}.\\
With this we obtain the following result for the worst-case error of QMC-integration on spaces of fractional smoothness:
Let $p\in(1,\infty)$ and $\alpha\in\big(\max\{1/p,1/2\},1\big)$. Then there exists a constant $c>0$ such that for any $N$-point cubature rule $Q$ we have
$$e^{\wor}(Q,\HH_{\alpha,s,p})\geq c N^{-\alpha} \ln(N)^{\frac{s-1}{2}}.$$
This lower bound was proved in \cite{UU15} for certain Bessel potential spaces; due to \eqref{idSpaces} it also holds for spaces of fractional smoothness. This result implies that
QMC-cubature rules using $(t,m,s)$-nets as discussed in \cite{GDMS24} are indeed optimal for spaces of fractional smoothness. 

Furthermore, we transfer the fractional Koksma-Hlawka inequalites from \cite{Dic08a} to the Bessel potential spaces of dominated mixed smoothness; more precisely, we obtain the following result.\\
Let $1<\alpha^{-1}<p<\infty$ and $1\leq q\leq \infty$. Then there exist constants $C,c>0$ such that for all pointsets $P=\{x_0,...,x_{N-1}\}\subseteq[0,1]^s$ we have
\begin{equation*}%\label{BessKH}
c D^\ast_{\alpha,s,p',q}(P) \le\sup_{ ||f||_{S^{\alpha}_{p}H([0,1]^s)}\leq1}\left|I_s(f) - Q_P(f) \right|  \le C D^\ast_{\alpha,s,p',q}(P),
\end{equation*}
where
\begin{equation*}
D^\ast_{\alpha,s,p',q}(P) := \left(\sum_{\emptyset \neq u \subseteq \{1,...,s\}}
 \left(\int_{[0,1]^{u}}
|\Delta_{\alpha}(\bst_u,u,P)|^{p'} \,\rd\bst_u \right)^{q/p'}
\right)^{1/q}
\end{equation*}
is the fractional discrepancy of the pointset $P$ with 
\begin{equation*}
\Delta_{\alpha}(\bst_u, u, P) := \alpha^{-|u|}
\prod_{j\in u}(1-t_j)^{\alpha} - \frac{1}{N} \sum^{N-1}_{n=0}
\prod_{j\in u}(x_{n,j}-t_j)_+^{\alpha-1}
\end{equation*}
being the fractional discrepancy function, see Corollary \ref{KHI1}.\\
Similarly to the Koksma-Hlawka inequalities proved in \cite{Tri10} for Besov spaces, we also get the following result. 
Let $1<\alpha^{-1}<p<\infty$ and $1\leq q\le \infty$. Then there exists a constant $C>0$ such that for all pointsets $P=\{x_0,...,x_{N-1}\}\subseteq[0,1]^s$ 
and any $f \in S^{\alpha}_{p}H([0,1]^s)$ we have
\begin{equation*}%\label{BessKH}
\left|I_s(f) - Q_P(f) \right|  \le C D^{BP}_{1-\alpha,s,p',q}(P) 
||f||_{S^{\alpha}_{p}H([0,1]^s)},
\end{equation*}
where
\begin{equation*}
D^{BP}_{1-\alpha,s,p',q}(P) := \left(\sum_{\emptyset \neq u \subseteq \{1,...,s\}}
 ||\Delta_1(\bst_u,u,P)||^{q}_{S^{1-\alpha}_{p'}H([0,1]^u)}
\right)^{1/q},
\end{equation*}
with 
$$\Delta_1(\bst_u, u, P) := 
\prod_{j\in u}(1-t_j) - \frac{1}{N} \sum^{N-1}_{n=0}
\prod_{j\in u} \mathbb 1_{[0, x_{n,j})}(t_j)
%(x_{n,j}-t_j)_+^{\alpha-1}
$$
being the classical local discrepancy function anchored in $\bsone$, see Corollary \ref{KHI2}.
%and 
%\begin{equation*}
%\Delta(\bst_u,u,P):=\prod_{j\in u}t_j-\frac 1N\sum_{n=0}^{N-1}\prod_{j\in u}\mathbb 1_{( x_{n,j},1]}(t_j)
%\end{equation*}
%is the discrepancy function.

In the Appendix we give an elaborated proof of a theorem by Triebel on function extensions for Bessel potential spaces \cite{Tri92}, tailored to our needs.

\subsection{Some Remarks on Notation}
Let us close the introduction with some remarks concerning notation. Throughout the paper we assume $1\leq p,q\leq \infty$. We write $\N$ for the natural numbers with $\N=\{1,2,3,...\}$ and for $s\in\N$ we put $S:=\{1,...,s\}$. For a set $A$ we denote by $|A|$ its cardinality if $A$ is finite and its Lebesgue measure if $A$ is an infinite measurable subset of $\R^s$. 
For a function $f:\dom (f)\to \C$ the support of $f$ is defined by 
$$\supp(f):=\overline{\{x\in\dom(f)|\, f(x)\neq 0\}}.$$
For a measurable subset $\Omega$ of $\R^s$ we denote by $L^p(\Omega)$ the space of all $p$-integrable functions with regard to the (restricted) Lebesgue measure; the norm in $L^p(\Omega)$ is denoted by $||\cdot||_{L^p(\Omega)}$. If the set $\Omega$ is obvious from the context, we only write $||\cdot||_{L^p}$.\\
Furthermore, if $A,B$ are normed spaces, we write $A\hookrightarrow B$ if $A$ is continuously embedded into $B$, i.e., if $A\subseteq B$ and the map $\iota:A\to B,\, x\mapsto x$ is continuous.\\
In inequalities, we use $C$ for generic constants that may vary between lines.

\section{The Unit Interval }

\subsection{Spaces of Fractional Smoothness}
Most of this section is based on \cite{GDMS24}.\\
%We endow $\Omega\subseteq \R$ with the usual Borel-$\sigma$-field, endowed with the Lebesgue measure.
%For $\Omega=[0,1]$ or $\Omega=\R_+$ and with Cauchy's formula for repeated integration in mind, we now define a fractional integral in the following way.
We endow $[0,1]\subseteq \R$ with the usual Borel-$\sigma$-field, endowed with the Lebesgue measure. With Cauchy's formula for repeated integration in mind, we now define a fractional integral in the following way.
\begin{definition}\label{fracInt}
For $\alpha\in (0,\infty)$ and a measurable function $ f:[0,1]\to\R$ we define 
$$\frac{d^{-\alpha }}{dx^{-\alpha}}f(x):=
\frac{1}{\Gamma(\alpha)}\int_0^xf(t)(x-t)^{\alpha-1} \, dt$$
for all $x\in[0,1]$, for which the integrand is Lebesgue-integrable.
\end{definition}
%For $p>\max\{1,\alpha^{-1}\}$ and $f\in L^p(\Omega)$ one gets $f(\cdot)(x-\cdot)^{\alpha-1}\in L^1([0,x])$ and therefore $\frac{d^{-\alpha }}{dx^{-\alpha}}f(x)$ is defined for all $x\in\Omega$.\\ %Hölders inequality for all $x\in (0,1]$.\\
%Applying a derivative to this formula yields the following definition of a fractional derivative.
We obtain fractional derivatives by first doing fractional integration and afterwards taking derivatives (of the right order).
%We start by defining the fractional derivative in the following way.
\begin{definition}\label{fracDer}
Let $\alpha \in[0,\infty)$, $\beta \in\N$ with $\beta-1\leq\alpha<\beta$ and 
$f:[0,1]\to\R$ measurable, such that $\frac{d^{\alpha-\beta }}{dx^{\alpha-\beta}}f(x)$
exists for all $x\in[0,1]$. If for an $x\in[0,1]$ the derivative 
$$\frac{d^\alpha }{dx^\alpha}f(x):= \frac{1}{\Gamma (\beta -\alpha)}\frac{d^\beta}{dx^\beta}\int_0^xf(y)(x-y)^{\beta-\alpha -1}\, dy$$
exists, then we call $\frac{d^\alpha }{dx^\alpha}f(x)$ the fractional derivative of f in $x$ in the sense of Riemann and Liouville.
\end{definition}
%\textbf{Convention:} We will use the fractional derivative and its extension by zero synonymously.\\
With these definitions in mind we will now introduce spaces of fractional smoothness.
\begin{definition}\label{fracSpace}
Let $1\leq\alpha^{-1}<p\leq\infty$. We put
$$\mathcal H_{\alpha ,p}:=\Big \{ f:[0,1]\to \R \mid \exists \tilde{f}\in L^p([0,1]) \ \exists \ c\in \R : f=c+\frac{1}{\Gamma(\alpha)}\int_0^1\tilde{f}(t)(\cdot -t)^{\alpha -1}_+\, dt \Big\}.$$
Here $(\cdot -t)^{\alpha -1}_+$ denotes the function, which is $(x -t)^{\alpha -1}$ for $x>t$ and $0$ elsewhere. 
\end{definition}
\begin{remark}\label{fracRem}
Due to Hölders inequality and $1\leq\alpha^{-1}<p\leq\infty$, we get for $\tilde f\in L^p([0,1])$ that $\tilde{f}(\cdot)(x-\cdot)^{\alpha -1}_+\in L^1([0,1])$ for all $x\in[0,1]$.
Moreover, for $\alpha^{-1}< p_1\leq p_2\leq\infty$ one gets $\mathcal H_{\alpha ,p_2}\subseteq \mathcal H_{\alpha ,p_1}$. % due to the corresponding statements for Lebesgue-spaces with finite measure.
\end{remark}
We now define for $\alpha>0$ the differential operator (cf. \cite{Dic08a})
$$D^\alpha f(t):=\frac{d^\alpha (f(t)-f(0))}{dt^\alpha}.$$
With this, we state some basic properties regarding the function spaces $\mathcal H_{\alpha ,p}$. More precisely, statement $(i)$ is from \cite{Dic08a}, statement $(ii)$ and $(iv)$ from \cite{GDMS24} and statement $(iii)$ from \cite{Mei25}.

\begin{lemma}\label{fracSpaceLemma}
Let $1\leq\alpha^{-1}<p\leq\infty$.
\begin{itemize}
\item[(i)] Each $f\in \mathcal H_{\alpha ,p}$ is Hölder continuous with Hölder exponent $\alpha-1/p$.
\item[(ii)] If $\tilde{f}\in L^p([0,1]),\ c\in\R$ and 
$$f:=c+\frac{1}{\Gamma(\alpha)}\int_0^1\tilde{f}(t)(\cdot -t)^{\alpha -1}_+\, dt,$$
then for almost all $x\in[0,1]$ the fractional derivative $D^\alpha f(x)$ exists and is equal to $\tilde{f}(x)$.\\
In particular, we have for all $f\in \mathcal H_{\alpha,p}$
$$f=f(0)+\frac{1}{\Gamma(\alpha)}\int_0^1D^\alpha f(t)(\cdot -t)^{\alpha -1}_+\, dt.$$
\item[(iii)] Let $f:[0,1]\to\R$ be absolutely continuous and $f'\in L^q([0,1])$ with $(1-\alpha)^{-1}<q\leq\infty$. 
Then $f\in \mathcal{H}_{\alpha,\infty}$.%, $D^\alpha f$ existiert fast überall und $D^\alpha f=\frac{d^{\alpha-1}}{dx^{\alpha-1}}f'$ fast überall.
\item[(iv)] For $0<\alpha<1$ there exist functions $f\in \mathcal H_{\alpha,p}$ whose fractional derivative $\frac{d^\alpha }{dx^\alpha}f$ is not in $L^p([0,1])$.
\end{itemize}
\end{lemma}
\begin{proof}
For $(i)$ we refer to \cite{Dic08a}, for $(ii)$ and $(iv)$ to \cite{GDMS24}.\\
$(iii)$ For $\alpha=1$ and $q=\infty$ the statement is trivial. Assume now $0<\alpha<1$.\\
One gets immediately
$$\frac{d^{\alpha-1}}{dx^{\alpha-1}}f'(x)=\frac{1}{\Gamma(1-\alpha)}\int_0^xf'(t)(x-t)^{-\alpha}\, dt\in\mathcal{H}_{1-\alpha,q}$$
and therefore due to $(i)$ the continuity of $\frac{d^{\alpha-1}}{dx^{\alpha-1}}f'$. This implies
$\frac{d^{\alpha-1}}{dx^{\alpha-1}}f'\in L^\infty([0,1])$, which leads to 
$$f(0)+\frac{1}{\Gamma (\alpha)}\int_0^x\big(\frac{d^{\alpha-1}}{dx^{\alpha-1}}f'(t)\big) (x-t)^{\alpha-1}\, dt\in\mathcal{H}_{\alpha,\infty}.$$
Using Fubini's theorem and Euler's integral representation of the beta function 
$$\int_0^1z^{\alpha-1}(1-z)^{-\alpha}\, dz=\Gamma(1-\alpha)\Gamma(\alpha),$$
one now gets
\begin{eqnarray*}
\begin{aligned}
&f(0)+\frac{1}{\Gamma (\alpha)}\int_0^x\big(\frac{d^{\alpha-1}}{dx^{\alpha-1}}f'(t)\big) (x-t)^{\alpha-1}\, dt\\
=&f(0)+\frac{1}{\Gamma (\alpha)\Gamma(1-\alpha)}\int_0^x\int_0^tf'(s)(t-s)^{-\alpha}\, ds\, (x-t)^{\alpha-1}\, dt\\
%=&f(0)+\frac{1}{\Gamma (\alpha)\Gamma(1-\alpha)}\int_0^x\int_s^xf'(s)(t-s)^{-\alpha} (x-t)^{\alpha-1}\, dt\, ds\\
=&f(0)+\frac{1}{\Gamma (\alpha)\Gamma(1-\alpha)}\int_0^xf'(s)\int_s^x(t-s)^{-\alpha} (x-t)^{\alpha-1}\, dt\, ds\\
=&f(0)+\frac{1}{\Gamma (\alpha)\Gamma(1-\alpha)}\int_0^xf'(s)\Gamma(\alpha)\Gamma(1-\alpha)\, ds\\
=&f(0)+\int_0^xf'(s)\, ds=f(x),
\end{aligned}
\end{eqnarray*}
where we used Lebesgue's version of the fundamental theorem of calculus. 
\end{proof}
We want to shortly mention, that the restriction $0<\alpha<1$ in $(iv)$ was not mentioned in the original paper \cite{GDMS24}, but is indeed needed.\\

We can now introduce a norm on $\mathcal H_{\alpha,p}$, which turns it into a Banach space.

\begin{definition}
Let $1\leq\alpha^{-1}<p\leq\infty$. We define the space of fractional smoothness 
$\mathcal H_{\alpha,p}$ as the vector space $\mathcal H_{\alpha,p}$ equipped with the norm 
$$||f||_{\alpha,p}:=|f(0)|+\Gamma(\alpha)^{-1}||D^\alpha f||_{L^p}.$$
\end{definition}
\begin{remark}
For $1\leq\alpha^{-1}<p\leq \infty$ fixed and $q\in(1,\infty]$ one can also define the equivalent norms 
$$||f||_{\alpha,p,q}:=\Big (|f(0)|^q+\Gamma(\alpha)^{-q}||D^\alpha f||^q_{L^p}\Big)^{1/q}$$
on $\mathcal H_{\alpha,p}$, with the usual modification for $q=\infty$.
\end{remark}

For the next definition cf., e.g., \cite{ZXZ09, LZZ22}.

\begin{definition}\label{RKBS}
Let $D\neq \emptyset$. A Banach space $B$ of functions $f:D\to \R$ is called a reproducing kernel Banach space (RKBS) if all the evaluation functionals
$\delta_x:B\to\R,\, x\in D$, are continuous, i.e. for each $x\in D$ there exists a constant $C_x>0$ with
$$|\delta_x(f)|=|f(x)|\leq C_x ||f||_B,\quad \text{for each } f\in B.$$
\end{definition}
%\fbox{Referenz angeben!}\\
It is an immediate consequence of this definition, that norm convergence in a reproducing kernel Banach space implies pointwise convergence.

\begin{proposition}\label{fracEmbLp}
Let $1\leq\alpha^{-1}<p\leq\infty$. Then $\mathcal{H}_{\alpha,p}$ is a reproducing kernel Banach space and
$$\mathcal{H}_{\alpha,p}\hookrightarrow L^p([0,1]).$$
\end{proposition}

\begin{proof}
We begin by showing that $\mathcal{H}_{\alpha,p}$ is a reproducing kernel Banach space. Let $f\in\mathcal H_{\alpha,p}$ and $x\in [0,1]$. 
With the help of Hölder's inequality we obtain for a suitable constant $C>0$
\begin{eqnarray*}
\begin{aligned}
&|f(x)|\leq|f(0)|+\frac{1}{\Gamma(\alpha)}\Big|\int_0^1D^\alpha f(t)(x-t)^{\alpha-1}_+\, dt\Big|\\
\leq&|f(0)|+\frac{1}{\Gamma(\alpha)}||D^\alpha(f)||_{L^p([0,1])}\cdot\Big(\int_0^1(x-t)^{(\alpha-1)p'}_+\, dt\Big)^{1/p'}\\
=&|f(0)|+\frac{1}{\Gamma(\alpha)}||D^\alpha(f)||_{L^p([0,1])}\cdot\Big(\frac{-(x-t)^{(\alpha-1)p'+1}}{(\alpha-1) p'+1}\big|_0^x\Big)^{1/p'}\\
%\leq&|f_n(0)-f(0)|+\frac{1}{\Gamma(\alpha)}||D^\alpha(f_n-f)||_{L^p([0,1])}\cdot\Big(\frac{x^{(\alpha-1)q+1}}{(\alpha-1) q+1}\Big)^{1/p'}\\
\leq&|f(0)|+\frac{1}{\Gamma(\alpha)}||D^\alpha(f)||_{L^p([0,1])}\cdot\Big(\frac{1}{(\alpha-1) p'+1}\Big)^{1/p'}\leq C ||f||_{\mathcal H_{\alpha,p}}.\\
\end{aligned}
\end{eqnarray*}	
In order to show the embedding $\mathcal{H}_{\alpha,p}\hookrightarrow L^p([0,1])$, observe that by Lemma \ref{fracSpaceLemma} $(i)$ 
one gets $\mathcal{H}_{\alpha,p}\subseteq C([0,1])\subseteq L^p([0,1])$. Thus the embedding is a consequence of the closed graph theorem.
\end{proof}
%\begin{remark}\label{ptwsRL01}
%We also showed that convergence in $\mathcal H_{\alpha,p}$ implies pointwise (and uniform) convergence.
%\end{remark}

\subsection{Multipliers}
%We will now recall multipliers and the Mihlin multiplier theorem.\\
For $f\in L^1(\R)$ and $\xi\in\R$ the Fourier transform of  $f$ is
$$(\mathcal{F}f)(\xi):=\frac{1}{\sqrt{2\pi}}\int_\R f(x)e^{-ix\xi}\, dx.$$
Furthermore, we denote the Schwartz space of rapidly decaying $C^\infty$-functions by $\mathcal S(\R)$ and the dual space of tempered distributions by $\mathcal S'(\R)$.\\
We mainly use \cite{Gra14} as a reference for the properties of the Fourier transform and Schwartz functions. 
\begin{definition}
Let $1\leq p<\infty$ and $m$ be a bounded, measurable function on $\R$. If the linear operator 
$$T_m(f):=\mathcal{F}^{-1}\Big ( m \cdot(\mathcal{F}f)\Big ),\ f\in\mathcal{S}(\R),$$
is bounded on $L^p(\R)$, then we call $m$ an $L^p$-multiplier (cf. \cite{Gra14}).
Here, the operator $T_m$ is first defined on the dense subset $\mathcal{S}(\R)$ and then extended onto $L^p(\R)$. We denote by $\mathcal{M}_p(\R)$ the space of all $L^p$-multipliers.
\end{definition}
%We now state the central theorem of this section.
\begin{theorem}[Mihlin multiplier theorem]\label{MMT}
Let $m$ be a bounded and differentiable function on $\R\setminus\{0\}$ such that there exists a
constant $C>0$ with
$$|\xi|\cdot\big|\frac{d}{d\xi}m(\xi)\big|\leq C$$
for all $\xi\in\R\setminus\{0\}$. Then $m\in\mathcal{M}_p(\R)$ for all $1<p<\infty$.
\end{theorem}
A proof of this theorem can  e.g. be found in \cite[Thm. 6.2.7]{Gra14} or in \cite[Thm. 4.23 \& Rem. 4.24]{Ab12}.\\
The following examples will play an important role in the proofs of our embedding theorems.

\begin{example}\label{m1m2}
For $\alpha\in (0,1)$ the functions $m_1,m_2:\R\setminus\{0\}\to \R$ with
$$m_1(\xi):=\frac{|\xi|^\alpha}{\big(1+\xi^2\big)^{\alpha/2}}\quad\text{and}\quad\ m_2(\xi):=\frac{\big(1+\xi^2\big)^{\alpha/2}}{1+|\xi|^\alpha}$$
are $L^p$-multipliers for all $1<p<\infty$\\
Indeed, clearly both functions are bounded. 
%The reason for this is that $m_1$ and $m_2$ are continuous on $\R$ and 
%$$\lim_{|\xi|\to\infty}m_1(\xi)=\lim_{|\xi|\to\infty}\frac{|\xi|^\alpha}{\big(1+\xi^2\big)^{\alpha/2}} =\lim_{|\xi|\to\infty}\frac{1}{\big(1/\xi^2+1\big)^{\alpha/2}}= 1 $$
%as well as
%$$\lim_{|\xi|\to\infty}m_2(\xi)=\lim_{|\xi|\to\infty} \frac{\big(1+\xi^2\big)^{\alpha/2}}{1+|\xi|^\alpha} =\lim_{|\xi|\to\infty}\frac{\big(1/\xi^2+1\big)^{\alpha/2}}{1/|\xi|^\alpha+1}= 1. $$
Furthermore, for all $\xi\neq 0$ we have
\begin{eqnarray*}
\begin{aligned}
|\xi|m_1'(\xi)=&|\xi|\frac{\sgn(\xi)\alpha|\xi|^{\alpha-1}\big(1+\xi^2\big)^{\alpha/2}-|\xi|^\alpha\alpha\xi\big(1+\xi^2\big)^{\alpha/2-1}}{\big(1+\xi^2\big)^{\alpha}}\\
=&\alpha \sgn(\xi)\frac{|\xi|^{\alpha}\big(1+\xi^2\big)^{\alpha/2}-|\xi|^{\alpha+2} \big(1+\xi^2\big)^{\alpha/2-1}}{\big(1+\xi^2\big)^{\alpha}}.\\
\end{aligned}
\end{eqnarray*}	
This function is continuous on $\R$ and both nominator and denominator behave for large $|\xi|$ like $|\xi|^{2\alpha}$, therefore $\lim_{|\xi|\to\infty}|\xi|\cdot |m_1'(\xi)|<\infty$ and the Mihlin multiplier theorem gives $m_1\in \mathcal{M}_p(\R)$.\\
The calculation for $m_2$ is analogously.\\
Now let $C_1,C_2\in\R$ and let $\phi:\R\setminus\{0\}\to\R$ be given by 
$$\phi(x):=\begin{cases}C_1, & \text{if } x>0,\\C_2, & \text{if }x<0. \end{cases}$$
Then the Mihlin multiplier theorem ensures that $\phi\cdot m_1,\,\phi\cdot m_2\in\mathcal{M}_p(\R)$ holds for $1<p<\infty$.
\end{example}
%\begin{remark}\label{m1m2r}
%Since constant factors do not influence the conditions of the Mihlin multiplier theorem we get for $C_1,C_2\in\R$ and $\phi:\R\setminus\{0\}\to\R$, 
%$$\phi(x):=\begin{cases}C_1, & \text{if } x>0,\\C_2, & \text{if }x<0 \end{cases}$$
%also $\phi\cdot m_1,\phi\cdot m_2\in\mathcal{M}_p(\R)$ for $1<p<\infty$.
%\end{remark}
%\fbox{Später hinter Mihlin setzen:}\\
\begin{comment}
    The following corollary of the Marcinkiewicz multiplier theorem can be found in \cite[Cor. 6.2.5]{Gra14}.
\begin{corollary}\label{MarMT}
Let $m$ be a bounded $s$-fold continuously differentiable function defined on $M:=\R^s\setminus\{\bsx\in\R^s|\,\exists i\in \{1,...,s\}: x_i=0\}$. Assume that for all $k\in \{1,...,s\}$ and all distinct $j_1,...,j_k\in \{1,...,s\}$ there exists a constant $C_{\{j_1,...,j_k\}}>0$ such that for all %$\xi_r\in\R\setminus\{0\}$ for $r\not\in\{j_1,...,j_k\}$ we have
$\xi\in M$ we have
$$|\partial_{j_1}...\partial_{j_k}m(\xi_1,...,\xi_n)|\leq C_{\{j_1,...,j_k\}}|\xi_{j_1}|^{-1}\cdot...\cdot |\xi_{j_k}|^{-1}.$$
Then $m$ in an $L^p$-multiplier for all $1<p<\infty$.
\end{corollary}
\end{comment}

\subsection{Bessel Potential Spaces}
We will now introduce and discuss a second type of function spaces, the Bessel potential spaces;
references are e.g. \cite{BL76} and \cite{Ab12}.
\begin{definition}
For $\alpha\in\R$ and $\xi\in\R$ put $<\xi>:=(1+\xi^2)^{1/2}$ and define 
\begin{eqnarray*}
\begin{aligned}
\mathcal{J}^\alpha:\mathcal{S}'(\R)&\to\mathcal{S}'(\R), \\
f&\mapsto \mathcal{F}^{-1}\Big(<\xi>^\alpha (\mathcal{F}f)(\xi)\Big).
\end{aligned}
\end{eqnarray*}	
\end{definition}
\begin{remark}\label{BessBij}
It is easily verified that for $\alpha_1,\alpha_2\in\R$ and $f\in\mathcal{S}'(\R)$ one gets
$$\mathcal{J}^{\alpha_1}\mathcal{J}^{\alpha_2}f=\mathcal{J}^{\alpha_1+\alpha_2}f.$$
Since $\mathcal{J}^0f=f$ for all $f\in\Sc'(\R)$, the mapping $\mathcal{J}^\alpha$ is a bijection on $\Sc'(\R)$ with inverse $\mathcal{J}^{-\alpha}$ for all $\alpha\in\R$.
\end{remark}
The Japanese bracket $<\cdot>$ has the following useful property.
\begin{lemma}\label{BessMult}
Let $\alpha>0$ and $1<p<\infty$. Then $<\cdot>^{-\alpha}\in\mathcal{S}'(\R)$ is an $L^p$-multiplier.
\end{lemma}
One can quickly prove this lemma by using the Mihlin multiplier theorem. \\
We can now define the second class of function spaces.

\begin{definition}
Let $\alpha\in\R$ and $1\leq p\leq\infty$. The Bessel potential space $H^\alpha_p(\R)$ is given by
$$H^\alpha_p(\R):=\big\{ f\in \mathcal{S}'(\R)| \mathcal{J}^\alpha f\in L^p(\R)\big\}$$
and is equipped with the norm
$$||f||_{H^\alpha_p(\R)}:=|| \mathcal{J}^\alpha f ||_{L^p}.$$ 
\end{definition}
Bessel potential spaces are Banach spaces and Remark \ref{BessBij} implies $H^0_p(\R)=L^p(\R)$.
We now state some properties of these spaces.
\begin{lemma}\label{BessEmbBess}
Let $\alpha_1,\alpha_2\in\R$ with $\alpha_1\leq\alpha_2$ and $1< p<\infty$. Then $H^{\alpha_2}_p(\R)\hookrightarrow H^{\alpha_1}_p(\R)$.
\end{lemma}
\begin{proof}
Let $f\in H^{\alpha_2}_p(\R)$. Observe that
$\mathcal J ^{\alpha_1}f=\mathcal{J}^{\alpha_1-\alpha_2}\mathcal{J}^{\alpha_2} f$
with $\mathcal{J}^{\alpha_2} f\in L^p(\R)$.\\
Since $<\cdot>^{\alpha_1-\alpha_2}$ is an $L^p$-multiplier due to Lemma \ref{BessMult}, there is a constant $C>0$ such that
$$||f||_{H^{\alpha_1}_p(\R)}=||\mathcal{J}^{\alpha_1} f ||_{L^p}=||\mathcal{J}^{\alpha_1-\alpha_2}\mathcal{J}^{\alpha_2} f ||_{L^p}\leq C ||\mathcal{J}^{\alpha_2} f ||_{L^p}=C ||f||_{H^{\alpha_2}_p(\R)}.$$
\end{proof}
\begin{corollary}\label{BessEmbLp}
Let $\alpha \geq 0$ and $1< p<\infty$. Then $H^\alpha_p(\R)\hookrightarrow L^p(\R).$
\end{corollary}
For the next statement and a short proof, see \cite[Thm. 2.3.3]{Tri83}.
\begin{lemma}\label{BessScDense}
Let $\alpha \in \R$ and $1\leq p<\infty$. Then $\mathcal{S}(\R)$ is a dense subset of $ H^\alpha_p(\R)$.
\end{lemma}
\begin{comment}
\begin{proof}
Since $\mathcal{J}^\alpha f\in L^p(\R)$ and $\mathcal{S}(\R)\subseteq L^p(\R)$ dense, there is a sequence $(g_n)_n\subseteq \mathcal{S}(\R)$ with $g_n\to \mathcal{J}^\alpha f$ in $L^p(\R)$. Remark \ref{BessBij} now implies
$ \mathcal{J}^{-\alpha}g_n \to f \text{ in } H^\alpha_p(\R)$.
Furthermore we have $(\mathcal{J}^{-\alpha}g_n)_n\subseteq \mathcal{S}(\R)$.\\
\end{proof}
\end{comment}

\begin{proposition}\label{BessCont}
Let $1<\alpha^{-1}<p<\infty$. Then $f\in H^\alpha_p(\R)$ is bounded and Hölder continuous with Hölder exponent $\alpha-1/p$.
\end{proposition}
%A proof of this can be found in [MAYM], as well as [???] (smth. by Triebel?!)
To prove Proposition \ref{BessCont}, use \cite[Sect. 2.3.2, Prop. 2 $(iii)$]{Tri83} and 
\cite[Sect. 2.7.1, Rem. 2]{Tri83}.\\
With the mentioned results from \cite{Tri83}, one also gets as an immediate consequence the following lemma.
\begin{lemma}\label{ptwsBess}
Let $1<\alpha^{-1}<p<\infty$. Then $H^\alpha_p(\R)$ is a reproducing kernel Banach space.
%Let $1<\alpha^{-1}<p<\infty$ and $f,f_n\in H^\alpha_p(\R)$ with $f_n\to f$. Then for all $x\in\R$
%$$f_n(x)\to f(x).$$
\end{lemma}
\begin{comment}
\begin{proof}
Let $f,f_n\in H^\alpha_p(\R)$ with $f_n\to f$ in $H^\alpha_p(\R)$. 
With $g_n:=\F^{-1}<\xi>^{\alpha}\F f_n$ and $g:=\F^{-1}<\xi>^{\alpha}\F f$
one has
$$|f_n(x)-f(x)|=|\big(\F^{-1}<\cdot>^{-\alpha}\F (g_n-g)\big)(x)|.$$
Since $f,f_n$ are continuous functions due to Proposition \ref{BessCont}, convolution with a mollifier
\linebreak \mbox{$(\psi_m)_m\subseteq C^\infty_c(\R)$} gives pointwise convergence. 
Therefore
$$\forall x\in\R\,\exists M\in\N:|f_n(x)-f(x)|\leq 2\cdot \Big|\big( \psi_M * \F^{-1}<\cdot>^{-\alpha}\F (g_n-g)\big)(x) \Big|.$$
With the convolution theorem and Hölder's inequality we now have
\begin{eqnarray*}
\begin{aligned}
&\Big|\Big( \psi_M * \F^{-1}\big(<\cdot>^{-\alpha}\F (g_n-g)\big)\Big)(x) \Big|\\
=&2\pi\Big| \F^{-1}\big(<\cdot>^{-\alpha}\cdot (\F\psi_M)\F (g_n-g)\big) (x)\Big|\\
=&2\pi\Big|\Big( \F^{-1}\big(<\cdot>^{-\alpha}\cdot (\F\psi_M)\big) *(g_n-g)\Big) (x)\Big|\\
\leq& ||\F^{-1}\big(<\cdot>^{-\alpha}\cdot (\F\psi_M)\big)||_{L^q} ||g_n-g||_{L^p}\to 0.
\end{aligned}
\end{eqnarray*}
\fbox{PROBLEM: M may depend on n, so the expression does not have to go to $0$!}
Better: Use embedding into Hölder space. Norm involves a sup. Norm. This implies pointwise convergence.\\
Check Triebel 1! Then this as an immediate consequence of hölder continuity.
\end{proof}    
\end{comment}
\begin{lemma}\label{BessNorm}
For $1<\alpha^{-1}<p<\infty$ we set for $f\in H^\alpha_p(\R)$
$$||f||_{\mathring H^{ \alpha}_p(\R)}:=|| \mathcal{F}^{-1}\Big(|\xi|^\alpha (\mathcal{F}f)(\xi)\Big)||_{L^p}.$$
Then
$$||\cdot||_{L^p}+||\cdot||_{\mathring H^{ \alpha}_p(\R)}$$
is an equivalent norm on $H^\alpha_p(\R)$.
\end{lemma}
A proof can be found in \cite[Thm. 6.3.2]{BL76}.\\
\begin{comment}
\begin{proof}
We use $m_1$ and $m_2$ from Example \ref{m1m2}.
Since $m_1\in \mathcal{M}_p(\R)$ one gets that $||\cdot||_{\mathring H^{ \alpha}_p(\R)}$ is well-defined because for $f\in H^{ \alpha}_p(\R)$ it implies
\begin{eqnarray*}
\begin{aligned}
||f||_{\mathring H^{ \alpha}_p(\R)}=&|| \mathcal{F}^{-1}\Big(|\xi|^\alpha (\mathcal{F}f)(\xi)\Big)||_{L^p}\\
=&|| \mathcal{F}^{-1}\Big(\frac{|\xi|^\alpha}{<\xi>^\alpha}\mathcal{F}\mathcal{F}^{-1}<\xi>^\alpha (\mathcal{F}f)(\xi)\Big)||_{L^p}\\
\leq&C||\mathcal{F}^{-1}<\xi>^\alpha (\mathcal{F}f)(\xi)||_{L^p}=C||f||_{H^{ \alpha}_p(\R)}.\\
\end{aligned}
\end{eqnarray*}	
Together with Corollary \ref{BessEmbLp} one gets the first inequality.\\
The second inequality follows in a similar fashion.
\begin{eqnarray*}
\begin{aligned}
|f||_{H^{ \alpha}_p(\R)}=&||\mathcal{F}^{-1}\Big(<\xi>^\alpha (\mathcal{F}f)(\xi)\Big)||_{L^p}\\
=&||\mathcal{F}^{-1}\Big(\frac{<\xi>^\alpha}{1+|\xi|^\alpha}\big(1+|\xi|^\alpha\big) (\mathcal{F}f)(\xi)\Big)||_{L^p}\\
%=&||\mathcal{F}^{-1}\Big(m_2\big(1+|\xi|^\alpha\big) (\mathcal{F}f)(\xi)\Big)||_{L^p}\\
=&||\mathcal{F}^{-1}\Big(m_2\mathcal{F}\mathcal{F}^{-1}\big(1+|\xi|^\alpha\big) (\mathcal{F}f)(\xi)\Big)||_{L^p}\\
\leq&C||\mathcal{F}^{-1}\big(\big(1+|\xi|^\alpha\big) (\mathcal{F}f)(\xi)\big)||_{L^p}\\
\leq&C\Big(||f||_{L^p}+||\mathcal{F}^{-1}\big(|\xi|^\alpha(\mathcal{F}f)(\xi)\big)||_{L^p} \Big)=C\Big(||f||_{L^p}+||f||_{\mathring H^{ \alpha}_p(\R)} \Big).
\end{aligned}
\end{eqnarray*}	
\end{proof}  
\end{comment}
The following result is due to Strichartz, 
cf. \cite[Thm. II2.1]{Str67}.
\begin{theorem}\label{BessMultip}
Let $1<\alpha^{-1}<p<\infty$ and $g\in H^\alpha_p(\R)$. Then there is a constant $C>0$ such that for all $f\in H^\alpha_p(\R)$ the inequality
$$||f\cdot g||_{H^\alpha_p(\R)}\leq C\cdot ||f||_{H^\alpha_p(\R)}$$
holds. We call a function with this property a \emph{pointwise multiplier} on $H^\alpha_p(\R)$.
\end{theorem}
We will now shortly introduce classical Sobolev spaces, since the connection between those spaces and Bessel potential spaces will be of importance later on.
\begin{definition}
For $n\in\N$ and $1\leq p<\infty$ we define
$$W^{n,p}(\R):=\{f\in L^p(\R)|\ \forall k\in\N,\ k\leq n:\ f^{(k)}\in L^p(\R)\}$$
and equip this space with the norm
$$||f||_{W^{n,p}}:=\Big(\sum_{k=0}^n\big|\big|f^{(k)}\big|\big|^p_{L^p}\Big)^{1/p};$$
the derivatives have to be understood in the distributional sense.
\end{definition}
\begin{theorem}\label{BessSob}
Let $n\in\N$ and $1< p<\infty$. Then
$$H^n_p(\R)=W^{n,p}(\R)$$
up to equivalence of norms.
\end{theorem}
A proof of this theorem can be found in \cite[Thm. 6.8]{Ab12}.

In order to get embedding results between fractional spaces in the sense of Riemann-Liouville and Bessel potential spaces, we will now define Bessel potential spaces on the unit interval.
\begin{definition}\label{Bess01}
For $\alpha \geq 0$ and $1\leq p<\infty$ we put
$$H^\alpha_p([0,1]):=
\big\{ f\in L^p([0,1])\,|\, \exists g\in H^\alpha_p(\R):g\mid _{[0,1]}=f \big\}$$
and equip this space with the norm
$$||f||_{H^\alpha_p([0,1])}:=\inf\big\{ ||g||_{H^\alpha_p} :\, g\mid_{[0,1]}=f\big\}.$$
\end{definition}
Bessel potential spaces on the unit interval are Banach spaces, see e.g. \cite[Prop. 3.2.3]{Tri83}. %(p. 193).\\
%We finish the section by showing some analogies of previous results for functions on the unit interval.
The following corollary is an immediate consequence of Lemma \ref{ptwsBess}.
\begin{corollary}\label{ptwsBess01}
%Let $1<\alpha^{-1}<p<\infty$ and $f,f_n\in H^\alpha_p([0,1])$ with $f_n\to f$ in $H^\alpha_p([0,1])$. Then for all $x\in \R$ holds
%$$f_n(x)\to f(x).$$
Let $1<\alpha^{-1}<p<\infty$. Then $H^\alpha_p([0,1])$ is a reproducing kernel Banach space.
\end{corollary}
%Corollary \ref{ptwsBess01} and the closed graph theorem imply the following corollary.
Corollary \ref{BessEmbLp} implies the following corollary.
\begin{corollary}\label{BessEmbLp01}
Let $\alpha\geq0$ and $1<p<\infty$. Then $H^\alpha_p([0,1])\hookrightarrow L^p([0,1]).$
\end{corollary}
%\begin{proof}
%$H^\alpha_p([0,1])\subseteq L^p([0,1])$ follows immediately from the embedding on $\R$, see corollary \ref{BessEmbLp}.\\
%Now let $f\in H^\alpha_p([0,1])$. Then there exists an extension $\tilde{f}\in H^\alpha_p(\R)$ of $f$ with 
%$$||\tilde{f}||_{H^\alpha_p(\R)}\leq 2||f||_{H^\alpha_p([0,1])}.$$
%Together with corollary \ref{BessEmbLp} one now gets
%$$||f||_{L^p([0,1])}\leq ||\tilde{f}||_{L^p(\R)}\leq C||\tilde{f}||_{H^\alpha_p(\R)}\leq 2C||f||_{H^\alpha_p([0,1])}.$$
%\end{proof}
%With similar arguments 
One can also prove embeddings between Bessel potential spaces on the unit interval, analogous to Lemma \ref{BessEmbBess}.

\subsection{Embeddings}
%We now want to investigate the relation between $\mathcal H_{\alpha ,p}$ and $H^\alpha_p([0,1])$.
%Therefore, we first have to introduce some new theory.\\
Our main goal in this section will be to prove for $1<\alpha^{-1}<p<\infty$ the embedding 
\begin{equation}\label{01emb}
H^\alpha_p([0,1])\hookrightarrow \mathcal H_{\alpha ,p },
\end{equation}
up to a few modifications we follow the approach of \cite{Mei25}.\\
%Therefore we first have to introduce some new function spaces.\\
We start by defining spaces $H^{\alpha }_{p,0}(\R)\subseteq H^\alpha _p(\R)$ and $H^{\alpha}_{p,0}([0,1])\subseteq H^\alpha_p([0,1])$ in which $f(0)=0$ holds (see Def. \ref{defRae}).
We will also introduce dense subsets $\Sc_0(\R)$ and $\Sc_0([0,1])$ of $H^{\alpha }_{p,0}(\R)$ and $H^{\alpha}_{p,0}([0,1])$, respectively, consisting of Schwartz functions (see Def. \ref{defRae} and Rem. \ref{RaeSubset}, as well as Thm. \ref{denseR} and Cor. \ref{dense01} for the density).
%We start by introducing the needed function spaces.
\begin{definition}\label{defRae}
For $1<\alpha^{-1}<p<\infty$ we put
\begin{eqnarray*}
\begin{aligned}
\Sc_0(\R)&:=\{f\in \Sc(\R) \mid f^{(n)}(0)=0\text{ for all }  n\in\N_0 \},\\
H^{\alpha }_{p,0}(\R)&:=\{f\in H^\alpha _p(\R)\mid f(0)=0\},\\
\Sc_0([0,1])&:=\{f\in C^\infty([0,1]) \mid \exists g\in \Sc_0(\R) : g|_{[0,1]}=f\},\\
H^{\alpha}_{p,0}([0,1])&:=\{f\in H^\alpha_p([0,1]) \mid \exists g\in H^{\alpha }_{p,0}(\R) : g|_{[0,1]}=f\}.
\end{aligned}
\end{eqnarray*}    
\end{definition}
Notice that due to Proposition \ref{BessCont} function evaluation is well-defined in $H_p^\alpha(\R)$ and thus also in $H_{p,0}^\alpha(\R)$ and $H^\alpha_{p,0}([0,1])$.
%The function values are to be read for the continuous representative of the equivalence classes.
\begin{remark}\label{RaeSubset}
For $1<\alpha^{-1}<p<\infty$ we have
$$
\Sc_0(\R)\subseteq H^{\alpha }_{p,0}(\R)\,\text{ and }\,
\Sc_0([0,1]) \subseteq H^{\alpha }_{p,0}([0,1]),
$$
due to the fact that the Schwartz functions are dense in all Bessel potential spaces, 
see Lemma~\ref{BessScDense}.
We can now equip all spaces from Definition \ref{defRae} with 
the corresponding $H^\alpha_p$-norm. Since convergence in $H^{\alpha }_{p}(\R)$ implies pointwise convergence (see Lemma \ref{ptwsBess}) and the space is complete, 
$H^{\alpha }_{p,0}(\R)$ is also a Banach space. The same holds for $H_p^\alpha([0,1])$.
\end{remark}
\begin{proposition}\label{denseR}
Let $1<\alpha^{-1}<p<\infty$. Then $\Sc_0(\R)$ lies dense in $ H^{\alpha }_{p,0}(\R) $.
\end{proposition}

\begin{proof}
Due to Lemma \ref{BessScDense} it is known that Schwartz functions lie dense in $H^{\alpha }_{p}(\R)$. We therefore can take for $f\in H^{\alpha }_{p,0}(\R) $ a sequence 
$(f_n)_n\subseteq \Sc(\R)$ with $f_n\to f$ in $H^{\alpha }_{p}(\R) $. Lemma \ref{ptwsBess} then implies $f_n(0)\to0$. We can now fix a $\psi\in C^\infty_0(\R)$ with $\psi(0)=1$ 
and define
$$g_n:=f_n-f_n(0)\psi.$$
With this $g_n\to f$ in $H^\alpha_p(\R)$ and $g_n(0)=0$ for all $n\in\N$. In order to get the vanishing derivatives in $x=0$, we modify these functions further. To do so, choose a 
$\rho\in C^\infty(\R)$ with $\rho(x)=0$ for $x\in [-1,1]$, $\rho(x)=1$ for $x\not\in [-2,2]$ and $\rho(x)\in [0,1]$ for all $x\in \R$. We can now define for $n\in\N$ and $\varepsilon>0$
$$h_{n,\varepsilon}(x):=\rho(x/\varepsilon)g_n(x).$$
By the Taylor expansion we now obtain for fixed $n\in\N$ the asymptotic bound $g_n(x)=O_n(x)$ as $x\to 0$. The difference $h_{n,\varepsilon}-g_n$ is supported in 
$[-2\varepsilon,2\varepsilon]$ and due to $|h_{n,\varepsilon}(x)|\leq |g_n(x)|$ for all $x\in\R$, we have
\begin{eqnarray*}
\begin{aligned}
&||h_{n,\varepsilon}-g_n||_{L^p(\R)}=||h_{n,\varepsilon}-g_n||_{L^p([-2\varepsilon,2\varepsilon])}\\
\leq&2||g_n||_{L^p([-2\varepsilon,2\varepsilon])}\leq 2\Big(\int_{-2\varepsilon}^{2\varepsilon} c_n|x|^p\,dx\Big)^{1/p}=\tilde c_n \varepsilon^{1+1/p}
\end{aligned}
\end{eqnarray*}
for some suitable constants $c_n,\tilde c_n,$ depending on $n$.
Hence 
$$||h_{n,\varepsilon}-g_n||_{L^p(\R)}=O_n(\varepsilon^{1+1/p}).$$
Next we examine the asymptotics of the derivative. Observe
$$(h_{n,\varepsilon}-g_n)'(x)= (\rho(x/\varepsilon)-1)g_n'(x)+\frac1\varepsilon \rho'(x/\varepsilon) g_n(x),$$
which is supported in $[-2\varepsilon,2\varepsilon]$.
%Due to $g_n(0)=0$ we have 
%$$g_n(x)=g_n'(0)x+O_n(x^2).$$
Substituting $x=\varepsilon t$ gives 
$$||(h_{n,\varepsilon}-g_n)'||_{L^p(\R)}^p=\varepsilon \int_{-2}^2\left|(\rho(t)-1)g_n'(t\varepsilon)+\rho'(t) \frac{g_n(t\varepsilon)}\varepsilon\right|^p\,dt.$$
For fixed $n\in\N$ and uniformly in $t\in[-2,2]$ we now have $g_n'(\varepsilon t)=g_n'(0)+O_n(\varepsilon)$ 
as well as $g_n(\varepsilon t)/\varepsilon=g_n'(0)t+O_n(\varepsilon)$ as $\varepsilon\to 0$. By defining
$$A_\rho(t):=\rho(t)-1+t\rho'(t),$$
we get for the integrand
$$(\rho(t)-1)g_n'(t\varepsilon)+\rho'(t) \frac{g_n(t\varepsilon)}\varepsilon=g_n'(0)A_\rho(t)+O_n(\varepsilon).$$
By the triangle inequality for the norm on $L^p([-2,2])$ one then gets
$$||(h_{n,\varepsilon}-g_n)'||_{L^p([-2,2])}=\varepsilon^{1/p} \left(|g_n'(0)| \, ||A_\rho||_{L^p([-2,2])}+O_n(\varepsilon) \right),$$
implying 
$$||(h_{n,\varepsilon}-g_n)'||_{L^p(\R)}=O_n(\varepsilon^{1/p}).$$
With Theorem \ref{BessSob} we now get $||h_{n,\varepsilon}-g_n||_{H^1_p(\R)} \to 0$ as $\varepsilon\to 0$ and the embedding 
$H^1_p(\R)\hookrightarrow H^\alpha_p(\R)$ implies
$||h_{n,\varepsilon}-g_n||_{H^\alpha_p(\R)} \to 0$ as $\varepsilon\to 0$.
These asymptotics are still dependent on $n$, choose therefore $\varepsilon_n>0$ such that 
the difference in the $H^\alpha_p(\R)$-norm is at most $1/n$. Now $h_{n,\varepsilon_n}\in\Sc_0(\R)$ since it vanishes in a neighborhood of zero and 
$$||h_{n,\varepsilon_n}-f||_{H^\alpha_p(\R)}\leq ||h_{n,\varepsilon_n}-g_n||_{H^\alpha_p(\R)}+||g_n-f||_{H^\alpha_p(\R)}\to 0,$$
if $n$ tends to infinity.
\end{proof}

\begin{corollary}\label{dense01}
Let $1<\alpha^{-1}<p<\infty$. Then $\Sc_0([0,1])$ lies dense in $ H^{\alpha }_{p,0}([0,1])$.
\end{corollary} 
\begin{definition}
For $\Omega\subseteq\R$ and $f:\Omega\to\R$ we define $\exto f:\R\to\R$ by
$$\exto f(x):=\begin{cases}  f(x),&x\in\Omega, \\    
0,&x\in\R\setminus\Omega.  \end{cases}\\$$
\end{definition}

\begin{lemma}\label{Lemma1}
Let $1<\alpha^{-1}<p<\infty$. Then $\exto(\cdot \,|_{\R_+})$ maps $H_{p,0}^\alpha(\R)$ into itself and
the mapping 
$$\exto(\cdot | _{\R_+}): H^\alpha_{p,0}(\R)\to H^\alpha_{p,0}(\R)$$
is continuous.
\end{lemma}
\begin{proof}
First consider $f\in \Sc_0(\R)$.\\
We have $\exto(f| _{\R_+})\in \Sc_0(\R)$ since $\exto(f| _{\R_+})$ is a $C^\infty$-function. Furthermore, the $L^p$-norms of $f$ and $f'$  under $\exto(\cdot| _{\R_+})$ do not get larger. Therefore we get $\exto(f| _{\R_+})\in H^1_p(\R)$ due to Theorem~\ref{BessSob}.\\
For $k\in\N$ we now set
\begin{eqnarray*}
\begin{aligned}
&a_0=1+10^{-k},\ &\lambda_0&=\frac{1}{2},\\
&a_1=-10^{-k}, \ &\lambda_1&=\frac{3}{2}\cdot 10^k+\frac{1}{2}.
\end{aligned}
\end{eqnarray*}
Since $a_0+a_1=1$ and $a_0\lambda_0+a_1\lambda_1=-1$ one gets with \cite[Thm. 4.5.2]{Tri92} the continuous linear operator %p. 223
\begin{eqnarray*}
\begin{aligned}
\exti: H^{\alpha}_p(\R_+)& \to  H^{\alpha}_p(\R),\\
\exti f(x)& = \begin{cases}  f(x),&x\in\R_+, \\    
\sum_{i=0}^1a_if(-\lambda_ix),&x\in\R_-  \end{cases}\\
\end{aligned}
\end{eqnarray*}
and therefore the inequality
$$||\exti (f| _{\R_+})||_{H^{\alpha}_p(\R)}\leq C||f|_{\R_+}||_{H^{\alpha}_p(\R_+)}$$
for some constant $C>0$ independent of $f$.
The space $H^\alpha_p(\R_+)$ is defined analogous to $H^\alpha_p([0,1])$, see Definition \ref{Bess01}.\\% (we will discuss these spaces in more detail later on).
For the convenience of the reader a version of the above mentioned theorem by Triebel tailored to our needs and a proof can be found in the Appendix.\\
%Our goal now is to get an inequality of the form
%$$||\exto(f\mid _{\R_+})||_{ H^{ \alpha}_p(\R)}\leq C||f||_{H^{\alpha}_p(\R)}.$$
Due to the definition of $\exti$ we now get
\begin{eqnarray*}
\begin{aligned}
&||\exto(f|_{\R_+})||_{\mathring H^{\alpha}_p(\R)}\\
\leq& ||\exti (f|_{\R_+})||_{\mathring H^{\alpha}_p(\R)}+|a_0|\cdot ||\exto(f|_{\R_+})(-\lambda_0 \ \cdot)||_{\mathring H^{\alpha}_p(\R)}\\
+&|a_1|\cdot ||\exto(f|_{\R_+})(-\lambda_1\ \cdot)||_{\mathring H^{ \alpha}_p(\R)}\\
=& ||\exti(f|_{\R_+})||_{\mathring H^{ \alpha}_p(\R)}+|a_0|\cdot \lambda_0^{ \alpha -1/p}\cdot||\exto(f|_{\R_+})||_{\mathring H^{ \alpha}_p(\R)}\\
+&|a_1|\cdot \lambda_1^{\alpha -1/p}\cdot||\exto(f|_{\R_+})||_{\mathring H^{\alpha}_p(\R)}.\\
\end{aligned}
\end{eqnarray*}
The last inequality for $g\in \Sc(\R)$ and $\lambda\in \R\setminus\{0\}$ is due to 
\begin{eqnarray*}
\begin{aligned}
&||g(\lambda\ \cdot)||_{\mathring H^{\alpha}_p(\R)}=||\F^{-1}\big(|\xi |^\alpha \F( g(\lambda\ \cdot))(\xi)\big) ||_{L^p} \\
=& ||\F^{-1}\big(|\xi |^\alpha \frac{1}{|\lambda|} (\F g)(\xi /\lambda)\big) ||_{L^p}=
|\lambda|^\alpha || \F^{-1}\big(\Big |\frac{\xi}{\lambda}\Big |^\alpha \frac{1}{|\lambda|} (\F g)(\xi /\lambda)\big) ||_{L^p}\\
=&|\lambda|^\alpha || (\F^{-1}(|\cdot |^\alpha (\F g)))(\lambda\ \cdot) ||_{L^p} = |\lambda|^{\alpha -1/p} || (\F^{-1}(|\cdot |^\alpha (\F g)))(\cdot) ||_{L^p}\\
=&|\lambda|^{\alpha -1/p} ||g||_{\mathring H^{\alpha}_p(\R)}.
\end{aligned}
\end{eqnarray*}
For $k\in \N$ sufficiently big and a fixed $\alpha\in(0,1)$ one now gets
$$|a_0|\cdot \lambda_0^{\alpha -1/p}+|a_1|\cdot \lambda_1^{\alpha -1/p}=:C_\alpha<1.$$
This is because $\alpha-1/p\in (0,1)$ and
\begin{eqnarray*}
\begin{aligned}
&|a_1|\cdot \lambda_1^{\alpha -1/p}=10^{-k}\cdot \Big(\frac{3}{2}\cdot 10^k+\frac{1}{2}\Big)^{\alpha-1/p}\\
\leq&10^{-k}\cdot \big(3\cdot 10^k\big)^{\alpha-1/p}\leq 3\cdot 10^{k(\alpha-1/p-1)}\overset{k\to\infty}{\to} 0.
\end{aligned}
\end{eqnarray*}
Moreover, since $\lambda_0^{\alpha -1/p}\in (0,1)$ we get $|a_0|\cdot \lambda_0^{\alpha -1/p}\in(0,1)$ for $k\in\N$ sufficiently big.\\
This leads to 
$$||\exto(f|_{\R_+})||_{\mathring H^{ \alpha}_p(\R)}-C_\alpha ||\exto(f|_{\R_+})||_{\mathring H^{ \alpha}_p(\R)} \leq||\exti(f|_{\R_+})||_{\mathring H^{ \alpha}_p(\R)}$$
and with this
$$||\exto(f|_{\R_+})||_{\mathring H^{ \alpha}_p(\R)}\leq\frac{1}{1-C_\alpha}||\exti(f|_{\R_+})||_{\mathring H^{ \alpha}_p(\R)}.$$
The definition of the norm on a subset gives us 
$$||f\mid _{\R_+}||_{H^{\alpha}_p(\R_+)}\leq ||f||_{H^{\alpha}_p(\R)}$$
and together with Lemma \ref{BessNorm} we now get
\begin{eqnarray*}
\begin{aligned}
||\exto(f|_{\R_+})||_{ H^{ \alpha}_p(\R)}\leq& C (||\exto(f|_{\R_+})||_{\mathring H^{ \alpha}_p(\R)}+||\exto(f|_{\R_+})||_{L^p})\\
\leq &C(||\exti(f|_{\R_+})||_{\mathring H^{ \alpha}_p(\R)}+||\exti(f|_{\R_+})||_{L^p})\\
\leq &C||\exti(f|_{\R_+})||_{H^{ \alpha}_p(\R)} \le C ||f|_{\R_+}||_{H^{\alpha}_p(\R_+)} \leq C  ||f||_{H^{\alpha}_p(\R)}.
\end{aligned}
\end{eqnarray*}
Take now $h\in H^{\alpha }_{p,0}(\R)$. We also assume $h\not\equiv 0$, since $0\in \Sc_0(\R)$.
Due to Proposition \ref{denseR} there now exists a sequence $(h_n)_n\subseteq  \Sc_0(\R)$ with $h_n \to h$ in $H^\alpha_p(\R)$ and due to the previously shown inequality $(\exto(h_n\mid _{\R_+}))_n$ is a Cauchy sequence in $H^\alpha_p(\R)$.\\
Since $H^\alpha_p(\R)$ is complete, $(\exto(h_n\mid _{\R_+}))_n$ must converge to some function in $H^\alpha_p(\R)$ and since convergence in $H^\alpha_p(\R)$ implies convergence in $L^p$  (see Prop. \ref{fracEmbLp}) $\exto(h\mid _{\R_+})$ must be the limit. \\
This implies $\exto(h\mid _{\R_+})\in H^\alpha_p(\R)$ and as an immediate consequence 
$\exto(h\mid _{\R_+})\in H^{\alpha }_{p,0}(\R)$.\\
Since 
$$||\exto(h_n\mid _{\R_+})||_{H^\alpha_p(\R)}\leq C||h_n||_{H^\alpha_p(\R)}$$
holds for all $n\in\N$, we get
\begin{eqnarray*}
\begin{aligned}
&||\exto(h\mid _{\R_+})||_{H^\alpha_p(\R)}\leq||\exto(h\mid _{\R_+})-\exto(h_n\mid _{\R_+})||_{H^\alpha_p(\R)}+||\exto(h_n\mid _{\R_+})||_{H^\alpha_p(\R)}\\
\leq &||h||_{H^\alpha_p(\R)}+C||h_n||_{H^\alpha_p(\R)}\leq C||h||_{H^\alpha_p(\R)}&
\end{aligned}
\end{eqnarray*}
for $n\in\N$ big enough.
%für $n\in\N$ groß genug (abhängig von $h$). Dabei nutzen wir die Konvergenz von $\exto(h_n\mid _{\R_+})$ und $h_n$ aus, sowie $h \not\equiv 0$ und die Stetigkeit der Norm. 
\end{proof}

\begin{lemma}\label{Lemma2}
For $1<\alpha^{-1}<p<\infty$ there exists a constant $C>0$ so that for all $f \in \Sc_0([0,1])$ there is an extension $\tilde f\in \Sc_0(\R)$ with $\supp (\tilde f)\subseteq\R_+$ and
$$||\tilde f||_{H^\alpha_p(\R)} \leq C || f||_{H^\alpha_p([0,1])}.$$
\end{lemma}
\begin{proof}
Let $f\in \Sc_0([0,1])$ with $f\not\equiv 0$ and $g \in H^\alpha_{p,0}(\R)$ be an extension of  $f$ with
$$||g||_{H^\alpha_p(\R)}\leq  \tilde C||f||_{H^\alpha_p([0,1])}$$
for some $\tilde C>0$ independent of $f$. 
We can assume due to Lemma \ref{Lemma1} that $\supp g\subseteq \R_+$. 
Theorem~\ref{BessMultip} applied on $h\in C^\infty_0(\R)$ with $\supp h\subseteq [-1,2]$ and $h(x)=1$ for $x\in [0,1]$ then allows us to assume $\supp g\subseteq [0,2]$.
Let now $f^s$ be an $\Sc(\R)$-extension of $f$ with $\supp f^s\subseteq[0,\infty)$. Such a function exists due to the definition of $\Sc_0([0,1])$ and Lemma \ref{Lemma1}. We now put 
$g^s:=g-f^s\in H^\alpha_p(\R)$ with $\supp(g^s)\subseteq [1,\infty)$.\\
For a function $\phi\in C^\infty_0(\R)$ with $\int_\R \phi(x)\,dx=1$, $\phi(x)\geq0$ for all $x\in\R$ and $\supp(\phi)\subset [-1,1]$, we call the sequence
of functions $(\phi_n)_n$ with $\phi_n(x):=n\phi(nx)$ a \emph{mollifier} or \emph{approximate identity}, see e.g. \cite{Gra14}.
Take now some mollifier $(\phi_n)_n\subseteq C_0^\infty(\R)$ with $\supp(\phi_n)\subseteq[0,1/n]$. \\
It holds that $(\phi_n~*~g^s)_n\subseteq~C^\infty(\R)$ and due to the bounded support of $g$ one gets $(\phi_n~*~g)_n\subseteq~C^\infty_0(\R)$, implying 
$(\phi_n~*~g^s)_n\subseteq~\mathcal S(\R)$.
We also have $\phi_n*g^s \to g^s$ in $H^\alpha_p(\R)$ due to convergence results for mollifiers and
the convolution theorem leading to
$$||\phi_n*g^s-g^s||_{H^\alpha_p(\R)}=||\phi_n * \F^{-1}<\xi>^{\alpha} \F g^s-\F^{-1}<\xi>^{\alpha} \F g^s||_{L^p}.$$
Furthermore $\supp(\phi_n)\subseteq[0,1/n]$ implies $\supp (\phi_n*g^s) \subseteq [1,\infty)$.\\
In particular, we get that the function $f^s+\phi_n*g^s=:f_n\in \Sc_0(\R)$ is a smooth extension of $f$ with $\supp(f_n)~\subseteq~[0,\infty)$ and 
$$f_n\to f^s+g^s=f^s+g-f^s=g\text{ in }H^\alpha_p(\R).$$
%Looking at the sequence $(f_n)_n\subseteq \Sc_0(\R)$ with $f_n:=f^s+\phi_n*g^s$, we get 
%$$f_n\to g \text{ in }H^\alpha_p(\R)$$
Therefore for $n_0\in\N$ big enough we have
$$||f_{n_0}||_{H^\alpha_p(\R)}\leq2||g||_{H^\alpha_p(\R)}\leq 2 \tilde C||f||_{H^\alpha_p([0,1])}.$$
Put now $\tilde{f}:=f_{n_0}$.\\
The case $f\equiv 0\in \Sc_0([0,1])$ is trivial, since $\tilde{f}\equiv 0\in \Sc_0(\R)$.
\end{proof}
\noindent
The idea now is to use Lemma \ref{Lemma2} to write the fractional derivative as a convolution. With this we afterwards prove the desired embedding \eqref{01emb} for
$ \Sc_0([0,1])$.
A density argument then implies $H^{\alpha }_{p,0}([0,1])\subseteq \mathcal H_{\alpha ,p }$. 
Since $H_p^{\alpha}([0,1])=H_{p,0}^\alpha([0,1]) \oplus \R$, one sees easily
$H^{\alpha }_{p}([0,1])\subseteq\mathcal H_{\alpha ,p }$. By using the closed graph theorem the proof is finished.
%We have now collected the tools needed in order to proof the main result of this chapter.
\begin{theorem}\label{embedding01}
Let $1<\alpha^{-1}<p<\infty$. Then
$$H^\alpha_p([0,1])\hookrightarrow \mathcal H_{\alpha ,p }.$$
\end{theorem}
\begin{proof}
The proof consists of three main steps.\\
\textbf{Step 1:} 
Show $\Sc_0([0,1])\hookrightarrow \mathcal H_{\alpha ,p }$.\\
We start by mentioning that due to Lemma \ref{fracSpaceLemma} $(iii)$, we have
$\Sc_0([0,1])\subseteq \mathcal{H}_{\alpha,p}$.\\
By Lemma \ref{Lemma2} we get an extension $\tilde{f} \in \Sc_0(\R)$ of $f\in \Sc_0([0,1])$ with $\supp(\tilde f)\subseteq [0,\infty)$ and
$$||\tilde f||_{H^\alpha_p(\R)} \leq C || f||_{H^\alpha_p([0,1])},$$
where $C>0$ is independent of $f$.
We get for $x\in[0,1]$
\begin{eqnarray*}
\begin{aligned}
D^\alpha f(x)=&\frac{1}{\Gamma (1-\alpha)}\frac{d}{dx}\int_0^x(f(t)-f(0))(x-t)^{-\alpha}\, dt\\
=&\frac{1}{\Gamma (1-\alpha)}\frac{d}{dx}\int_0^xf(t)(x-t)^{-\alpha}\, dt
%=&\frac{1}{\Gamma (1-\alpha)}\frac{d}{dx}\int_0^1 f(t)(x-t)_+^{-\alpha}\, dt\\
=\frac{1}{\Gamma (1-\alpha)}\frac{d}{dx}\int_{-\infty}^\infty \tilde{f}(t)(x-t)_+^{-\alpha}\, dt\\
=&\frac{1}{\Gamma (1-\alpha)}\frac{d}{dx} (\tilde{f}*(\cdot)_+^{-\alpha})(x)
%=&\frac{1}{\Gamma (1-\alpha)}\F^{-1}\F\frac{d}{dx} (\tilde{f}*(\cdot)_+^{-\alpha})(x)\\
=\frac{1}{\Gamma (1-\alpha)}\F^{-1}\F \Big(\frac{d\tilde{f}}{dx}*(\cdot)_+^{-\alpha}\Big)(x).
\end{aligned}
\end{eqnarray*}
Here we use in the last equality $(\cdot)_+^{-\alpha}\in\Sc'(\R)$ which can quickly be shown, as well as the well known rule on the derivative of the convolution of a 
Schwartz function and a tempered distribution, see e.g. \cite{Gra14}.
By the convolution theorem it follows that
\begin{eqnarray*}
\begin{aligned}
&\frac{1}{\Gamma (1-\alpha)}\F^{-1}\F\Big (\frac{d\tilde{f}}{dx}*(\cdot)_+^{-\alpha}\Big)(x)
=\frac{1}{\Gamma (1-\alpha)}\F^{-1}\Big (\sqrt{2\pi}\F\Big(\frac{d\tilde{f}}{dx}\Big)(\xi)\cdot\F((\cdot)_+^{-\alpha})(\xi)\Big)(x)\\
=&\frac{1}{\Gamma (1-\alpha)}\F^{-1} \big(i\xi \sqrt{2\pi} (\F\tilde f)(\xi)\cdot\F((\cdot)_+^{-\alpha})(\xi)\big)(x).
\end{aligned}
\end{eqnarray*}
From \cite[p. 400]{BN71} we get
$$\F((\cdot)_+^{-\alpha})(\xi)=\frac{\Gamma(1-\alpha)}{\sqrt{2\pi}}|\xi|^{\alpha -1}e^{-i\pi \frac{1-\alpha}{2}\sgn(\xi)},$$
leading to
\begin{eqnarray*}
\begin{aligned}
D^\alpha f=&\frac{1}{\Gamma (1-\alpha)}\F^{-1} \big(i\xi \sqrt{2\pi} (\F\tilde f)\F((\cdot)_+^{-\alpha})\big)\\
=&i\F^{-1}\big (\xi (\F \tilde{f}) 
|\xi|^{\alpha -1}e^{-i\pi \frac{1-\alpha}{2}\sgn(\xi)}\big )\\
=&i\F^{-1}\big (|\xi|^\alpha (\F \tilde{f}) 
\sgn(\xi)e^{-i\pi \frac{1-\alpha}{2}\sgn(\xi)}\big )\\				
=&i\F^{-1}\big (\frac{|\xi|^\alpha}{<\xi>^{\alpha}}\sgn(\xi)e^{-i\pi \frac{1-\alpha}{2}\sgn(\xi)}\F\F^{-1}<\xi>^{\alpha} (\F \tilde{f}) \big ).
\end{aligned}
\end{eqnarray*}
The required norm-inequality now follows by the Mihlin multiplier theorem. Example \ref{m1m2} implies $\frac{|\xi|^\alpha}{<\xi>^{\alpha}}\in \mathcal{M}_p(\R)$ and that the factor $\sgn(\xi)e^{-i\pi \frac{1-\alpha}{2}\sgn(\xi)}$ can be neglected; therefore
$$\frac{|\xi|^\alpha}{<\xi>^{\alpha}}\sgn(\xi)e^{-i\pi \frac{1-\alpha}{2}\sgn(\xi)} \in \mathcal{M}_p(\R).$$
This implies the existence of a $\widetilde{C}>0$ such that
$$||D^\alpha f ||_{L^p}
\leq\widetilde{C} \big| \big|i\F^{-1}<\xi>^{\alpha} \F (\tilde{f}) \big ) \big| \big|_{L^p}=\widetilde{C}||\tilde{f}||_{H^{ \alpha}_p(\R)}.
$$
Altogether we obtain for $f\in \mathcal S_0([0,1])$
\begin{eqnarray*}
\begin{aligned}
&||f||_{\mathcal H_{\alpha ,p }}=|f(0)|+\frac{1}{\Gamma (\alpha)}||D^\alpha f||_{L^p}=\frac{1}{\Gamma (\alpha)}||D^\alpha f||_{L^p}\\
\leq& \frac{\widetilde{C}}{\Gamma (\alpha)}||\tilde{f}||_{H^{ \alpha}_p(\R)} \leq  \frac{C\cdot\widetilde{C}}{\Gamma (\alpha)} ||f||_{H^{ \alpha}_p([0,1])}.
\end{aligned}
\end{eqnarray*}
\textbf{Step 2:} Show $H^{\alpha }_{p}([0,1])\subseteq \mathcal H_{\alpha ,p }$.\\
We use a density argument.
Take $f\in H^{\alpha }_{p,0}([0,1])$ and $(f_n)_n\subseteq \Sc_0([0,1])$ a sequence converging to $f$ which exists due to Corollary \ref{dense01}.
Now $(f_n)_n$ is a Cauchy sequence in $H^{\alpha }_{p}([0,1])$ and due to Step 1 it is also a Cauchy sequence in $\mathcal H_{\alpha ,p }$. Therefore a limit $g\in \mathcal H_{\alpha ,p }$ must exist.\\
With Proposition \ref{fracEmbLp} and Corollary \ref{ptwsBess01} we now get for all $x\in [0,1]$
$$|f(x)-g(x)|\leq |f(x)-f_n(x)|+|f_n(x)-g(x)|\to 0,$$	
hence $f=g\in\mathcal H_{\alpha,p}$.
One can now conclude $H^{\alpha}_p([0,1])\subseteq \mathcal H_{\alpha ,p }$.\\
For this, use that constant functions on $[0,1]$ are in $H^{\alpha}_p([0,1])$, since for all $c\in\R$ there exists a function $f_c\in C^\infty_0(\R)\subseteq H^{\alpha}_p(\R)$ with $f_c(x)=c$ for all $x\in[0,1]$.\\
We now get for  $f\in H^{\alpha}_p([0,1])$
$$f=g+c,\ g\in H^{\alpha }_{p,0}([0,1]),\ c=f(0)\in\R.$$
This implies
$$f=f(0)+ \frac{1}{\Gamma(\alpha)}\int_0^1D^\alpha g(t)(\cdot -t)^{\alpha -1}_+dt \in\mathcal{H}_{\alpha,p}.$$
\textbf{Step 3:} Show $H^{\alpha }_{p}([0,1])\hookrightarrow \mathcal H_{\alpha ,p }$.\\
We use the closed graph theorem.
Take $(f_n)_n\subseteq H^{\alpha}_p([0,1])$ a sequence converging in $H^{\alpha}_p([0,1])$ to $f\in H^{\alpha}_p([0,1])$ and in $\mathcal H_{\alpha ,p }$ to $g\in \mathcal H_{\alpha ,p }$.
We get for all $x\in[0,1]$ as before
$$|f(x)-g(x)|\leq |f(x)-f_n(x)|+|f_n(x)-g(x)|\to 0.$$
The closed graph theorem now implies the continuity of the identical embedding of 
$H^{\alpha}_p([0,1])$ into $ \mathcal H_{\alpha ,p }$.
\end{proof}
From \cite{GDMS24} the following embedding is already known.

\begin{theorem}$\cite[Thm.~7]{GDMS24}$. 
Let $1<\alpha^{-1}<p<\infty$. Then
$$\mathcal{H}_{\alpha,p}\hookrightarrow H^\alpha_p([0,1]).$$
\end{theorem}
A proof can be found in the paper. 
We therefore get the following corollary.
\begin{corollary}\label{equal01}
Let $1<\alpha^{-1}<p<\infty$. Then
$$H^\alpha_p([0,1])= \mathcal H_{\alpha ,p }$$
up to equivalence of norms.
\end{corollary}

\begin{comment}
Therefore most properties can be transferred between the spaces. With Remark \ref{fracRem} and Lemma \ref{BessEmbBess} we get, e.g. 
\begin{corollary}\label{C1}
Let $1<\alpha^{-1}< p_1\leq p_2 <\infty$. Then
$$H^\alpha_{p_2}([0,1])\subseteq H^\alpha_{p_1}([0,1]).$$
\end{corollary}	
\begin{corollary}\label{C2}
Let $\alpha_1,\alpha_2\in (0,1),\ \alpha_1\leq\alpha_2$ and $\alpha_1^{-1}<p<\infty$. Then
$$\mathcal{H}_{\alpha_2,p}\subseteq \mathcal{H}_{\alpha_1,p}.$$
\end{corollary}    
\end{comment}

\section{The Unit Cube}

\subsection{Spaces of Fractional Smoothness}

Let $s\in \NN$, $\alpha>0$, and $p>\alpha^{-1}$.
We define the function space $\HH_{\alpha,s,p}$ as the vector space
\begin{equation}\label{def_frac_H_alpha_s_p}
\HH_{\alpha,s,p} := \{ f:[0,1]^s \to \RR \,|\,
\forall u \subseteq S \,\exists \tilde{f}_u \in
L_p([0,1]^{u}) : f =
\Phi((\tilde{f}_u)_{u\subseteq S}) \},
\end{equation}
where
\begin{equation}
\label{Phi}
\Phi((\tilde{f}_u)_{u\subseteq S})(\bsx) :=
\sum_{u\subseteq S} \Gamma(\alpha)^{-|u|}
\int_{[0,1]^{u}} \tilde{f}_u(\bst_u) \prod_{j\in u}
(x_j-t_j)^{\alpha-1}_+ \, \rd \bst_u.
\end{equation}
Here we used the convention that $L_p([0,1]^{\emptyset}) :=\RR$
and
\begin{equation*}
\int_{[0,1]^{\emptyset}} \tilde{f}_\emptyset(\bst_\emptyset)
\prod_{j\in\emptyset} (x_j-t_j)^{\alpha-1}_+\,\rd \bst_\emptyset
:= \tilde{f}_\emptyset \in \RR\,;
\end{equation*}
later the norm 
$\|\cdot\|_{L_p([0,1]^{\emptyset})}$ should be interpreted as the absolute value on $\RR$.
Note that for $\alpha^{-1} < p_1 \le p_2 \le \infty$ we have $\HH_{\alpha,s,p_2} \subseteq \HH_{\alpha,s,p_1}$.

For $\bsx, \bsy\in [0,1]^s$ we write 
$(\bsx_u,\bsy)$ for $(z_1,\ldots,z_s)$
if $z_j=x_j$ for all $j\in u$ and $z_j = y_j$ otherwise.
For a function $f:[0,1]^s\to \RR$, $u\subseteq S$, and 
$\bst_u\in [0,1]^{u}$ we put
\begin{equation*}
f_{\forksnot, u}(\bst_u) := \sum_{v\subseteq u} (-1)^{|u\setminus v|}\,
f(\bst_v, \bszero).
\end{equation*}
The term $f_{\forksnot, u}$ is the $u$-th term of the
\emph{anchored decomposition} (also known as \emph{cut-HDMR}, cf. \cite{RA1999}) 
$$
f=\sum_{u\subseteq S} f_{\forksnot, u}
$$ 
of
$f$ with anchor in $\bszero$, i.e., $f_{\forksnot, u}$ depends only on the variables $x_j$, $j\in u$, and  $f_{\forksnot, u}(\bsx) = 0$ if there exists  some $j\in u$ with $x_j=0$, see \cite{KSWW10b}.
We define now for $u\subseteq S$ (formally) the fractional
derivative
\begin{equation}
\label{partfracder}
D^{(\bsalpha_u, \bszero_{S\setminus u})}f(\bst_u, \bszero) :=
\prod_{j\in u} \frac{\partial^{\alpha}}{\partial t^\alpha_j}
f_{\forksnot, u}(\bst_u),
\end{equation}
where
$\frac{\partial^{\alpha}}{\partial t_j^\alpha}$ is the $j$-th partial
derivative corresponding to the fractional derivative of order
$\alpha>0$ in Definition~\ref{fracDer}.
Here we use the convention
\begin{equation*}
D^{(\bsalpha_\emptyset, \bszero_{S})}f(\bst_\emptyset,\bszero) 
:= f(\bszero).
\end{equation*}

Any order of applying
the partial fractional derivatives in (\ref{partfracder}) to a
function $f\in \HH_{\alpha,s,p}$ leads to the same result,
see, e.g., the proof of \cite[Lemma~9(ii)]{GDMS24}.
Hence an
analogon of the classical theorem of Schwarz for mixed classical
partial derivatives holds also for mixed fractional partial derivatives.

Statement~(i) of the next lemma was verified in \cite[Prop.~1]{Dic08a}
and statement~(ii) in \cite[Lemma~8]{GDMS24}.

\begin{lemma}\label{s-Derivative}
Let $s\in \NN$, $\alpha\in (0,1]$, and $p>\alpha^{-1}$.
\begin{itemize}
\item[{\rm (}i{\rm )}]
Each $f\in \HH_{\alpha,s,p}$ is H\"older continuous with H\"older
exponent $\alpha-1/p$.
\item[{\rm (}ii{\rm )}]
If $f\in \HH_{\alpha,s,p}$ and
$f = \Phi((\tilde{f}_u)_{u\subseteq S})$, where
$\tilde{f}_u \in L_p([0,1]^{u})$ for all $u\subseteq S$,
then for each $u \subseteq S$ and almost all
$t_u\in [0,1]^{u}$ the fractional derivative
$D^{(\bsalpha_u, \bszero_{S\setminus u})}f(\bst_u, \bszero)$
exists and is equal to $\tilde{f}_u(t_u)$.
In particular, we have for all $f\in \HH_{\alpha,s,p}$ that
\begin{equation*}
f = \Phi((D^{(\bsalpha_u, \bszero_{S\setminus u})}
f(\cdot, \bszero))_{u\subseteq S}).
\end{equation*}
Moreover, for each $u\subseteq S$ the $u$-th term of the anchored decomposition
of $f$ with anchor in $\bszero$ is given by
\begin{equation}\label{anchored_decomp}
f_{\forksnot, u} = \Gamma(\alpha)^{-|u|} \int_{[0,1]^u}  
D^{(\bsalpha_u, \bszero_{S\setminus u})} f(\bst_u,  \bszero) \prod_{j\in u} (x_j-t_j)_+^{\alpha -1}
\, {\rm d} \bst_u,
\end{equation}
i.e., the right hand side of \eqref{Phi} constitutes the anchored decomposition of $f$.
\end{itemize}
\end{lemma}

\begin{definition}\label{Def:Space_Frac_Smooth_Mult}
For $p,q\in [1,\infty]$ and $\alpha\in (1/p,1]$  let the \emph{space of fractional smoothness}
$\HH_{\alpha,s,p,q}$ be the vector space $\HH_{\alpha,s,p}$ defined as in \eqref{def_frac_H_alpha_s_p},
endowed with the
norm
\begin{equation*}
\|f\|_{\alpha, s,p,q} := \left( \sum_{u\subseteq S}
\Gamma(\alpha)^{-q|u|} \|\tilde{f}_u\|^{q}_{L_p}
\right)^{1/q}
\end{equation*}
\end{definition}

Note that the norm is well-defined due to Lemma~\ref{s-Derivative} (ii). 
It is easily seen that $\HH_{\alpha,s,p,q}$
is a complete normed space. 
%and actually it is even a reproducing kernel Banach space
%in the sense of \cite{GSGP14, LZZ19}, see \cite{GHLW20}.

\begin{remark}\label{s:q_egal}
Note that different choices of $q$ lead to equivalent norms
on the vector space $\HH_{\alpha,s,p}$.  Therefore we often suppress the reference to $q$ and write simply $\HH_{\alpha,s,p}$  instead of $\HH_{\alpha,s,p,q}$.
\end{remark}

%&&&&&&&&&&&&&&&&&&&&&&&&&&&&&&&&&&&&&&&&&&&&&&&&&&&&&&&&&&&&&&&&
%&&&&&&&&&&&&&&&&&&&&&&&&&&&&&&&&&&&&&&&&&&&&&&&&&&&&&&&&&&&&&&&&&&

\subsection{Multipliers}
%We will now recall multipliers and the Mihlin multiplier theorem.\\
For $f\in L^1(\R^s)$ and $\bsxi\in\R^s$ the Fourier transform of  $f$ is
$$(\mathcal{F}f)(\bsxi):=\frac{1}{(2\pi)^{s/2}}\int_{\R^s}f(\bsx)e^{-i\bsx\cdot\bsxi}\, d\bsx,$$
see \cite{Gra14}. Similarly as before, we denote the Schwartz space of rapidly decaying $C^\infty$-functions by $\mathcal S(\R^s)$ and the dual space of tempered distributions by $\mathcal S'(\R^s)$.\\
%We mainly use \cite{Gra14} as a reference for the properties of the Fourier transform and Schwartz functions. 
\begin{definition}
Let $1\leq p<\infty$ and $m$ be a bounded, measurable function on $\R^s$. If the linear operator 
$$T_m(f):=\mathcal{F}^{-1}\Big ( m \cdot(\mathcal{F}f)\Big ),\ f\in\mathcal{S}(\R^s),$$
is bounded on $L^p(\R^s)$, then we call $m$ an $L^p$-multiplier, cf. \cite{Gra14}.
Here, the operator $T_m$ is first defined on the dense subset $\mathcal{S}(\R^s)$ and then extended onto $L^p(\R^s)$. We denote by $\mathcal{M}_p(\R^s)$ the space of all $L^p$-multipliers.
\end{definition}

The following corollary of the Marcinkiewicz multiplier theorem can be found in \cite[Cor.~6.2.5]{Gra14}.
\begin{corollary}\label{MarMT}
Let $m$ be a bounded $s$-fold continuously differentiable function defined on $M~:=\R^s\setminus\{\bsx\in\R^s|\,\exists i\in \{1,...,s\}: x_i=0\}$. Assume that for all $k\in \{1,...,s\}$ and all distinct $j_1,...,j_k\in \{1,...,s\}$ there exists a constant $C_{\{j_1,...,j_k\}}>0$ such that for all %$\xi_r\in\R\setminus\{0\}$ for $r\not\in\{j_1,...,j_k\}$ we have
$\bsxi\in M$ we have
$$|\xi_{j_1}|\cdot...\cdot |\xi_{j_k}|\cdot |\partial_{j_1}...\partial_{j_k}m(\xi_1,...,\xi_n)|\leq C_{\{j_1,...,j_k\}}.$$
Then $m$ is an $L^p$-multiplier for all $1<p<\infty$.
\end{corollary}
We will now shortly discuss an important example of a multiplier.

\begin{example}\label{ex2}
    Let $\emptyset \neq u\subseteq S$, $\alpha\in (0,1)$ and define $M_u:=\R^u\setminus\{\bsx\in\R^u|\,\exists i\in \{1,...,u\}~:~x_i=~0~\}$. We define
    $$m_0:M_u\to \R,\quad \bsxi \mapsto \prod_{j\in u}\sign(\xi_j)e^{-i\pi\frac\alpha 2 \sign(\xi_j)}$$
    and for all $j\in u$
    $$m_j:M_u\to \R,\quad \bsxi \mapsto \frac{|\xi_j|^{1-\alpha}}{<\xi_j>^{1-\alpha}}.$$
    By Corollary \ref{MarMT} it is easily seen, that these are $L^p$-multipliers on $\R^u$ for $1<p<\infty$. Therefore this also holds for $m:=\prod_{j\in u\cup\{0\}}m_j$.
\end{example}

\subsection{Bessel Potential Spaces}

Let $1\leq p\leq\infty$ and $\alpha \in \real$.
The \emph{Bessel potential space} or \emph{(fractional) Sobolev space with dominating mixed smoothness}
$S^\alpha_pH(\real^s)$ is  the collection of all tempered distributions $f
\in \mathcal{S}'(\real^s)$ such that
\[
\| \, f \, \|_{S^{\alpha}_p H(\real^s)}: = \Big\|\, \cfi
\Big[\prod\limits_{i=1}^s (1+|\xi_i|^2)^{\alpha/2}\, \cf f (\xi)
\Big](\, \cdot \, )\Big\|_{L_p(\real^s)}
\]
is finite.  Similarly as in the univariate case, see Corollary \ref{BessEmbLp}, we have $S^\alpha_pH(\real^s) \subseteq L_p(\real^s)$ for $\alpha \ge 0$, $1<p<\infty$.

If we take $\alpha =1$, then we obtain $S^{1}_p H(\real^s) = S^{1}_p W(\real^s)$ in the sense of equivalent norms,
where $S^{1}_p W(\real^s)$ denotes the first order Sobolev space of dominating mixed smoothness.
This space is the collection of all $f\in L_p (\real^s)$ such that
\[
 \|f\|_{S^{1}_p W(\real^s)} := \|f\|_{L_p(\real^s)} + \sum_{{{\bf m}\in \N_0^s} \atop \|{\bf m}\|_\infty \le 1}
 %0 \le m_i \le 1~ \forall i
 \|D^{\bf m} f\|_{L_p (\real^s)}
<\infty\, .
 \]
Here $D^{\bf m} f$ is understood as the corresponding distributional derivate of $f$.
In many problems the mixed derivative $D^{(1, 1, \ldots \, 1)} f$ is the dominating one, which is  the reason for the name.

Since we want to deal with spaces on the unit cube $[0,1]^{s}$, 
we shall apply restrictions  to switch from $\real^{s}$  to the cube.  
Let $\alpha \ge 0$. 
The space
$S^{\alpha}_{p}H([0,1]^s)$ is the space of all $f\in
   L_p([0,1]^s)$ such that there exists a $g\in
   S^{\alpha}_{p}H(\real^s)$ satisfying $f = g|_{[0,1]^s}$. It is endowed  with the quotient norm
   $$
 \|\,f\,\|_{S^{\alpha}_{p}H([0,1]^s)} := \, \inf\{\|\, g \,\|_{S^{\alpha}_{p}H(\real^s)} \,|\, g\in S^{\alpha}_{p}H(\real^s)\,,\, g_{|_{[0,1]^s}} =f\}\,.
   $$

Next we restate \cite[Lemma~11]{GDMS24}.

\begin{lemma}\label{Win_Lemma_13}
 Let $s\ge 1$,  $\alpha \ge 0$, and $1 < p<\infty$. Then the following statements hold. 
 \\
{\rm (i)} 
%Let $\alpha \in [0,1]$ and $1 < p<\infty$.Then 
$ S^{\alpha}_{p}H([0,1]^s)$ is a Banach space.
\\
{\rm (ii)} The scale $ S^{\alpha}_{p}H([0,1]^s)$
is monotone in $p$ and $\alpha$, i.e., for any $\varepsilon >0$
we have 
\[
 S^{\alpha+\varepsilon}_{p}H([0,1]^s) \hookrightarrow  S^{\alpha}_{p}H([0,1]^s)
\]
 and 
 \[
 S^{\alpha}_{p+\varepsilon} H([0,1]^s) \hookrightarrow  S^{\alpha}_{p}H([0,1]^s)\, .
\]
 \end{lemma}

\subsection{Tensorization}

For tensor products of general Banach spaces we refer to \cite{DF93} and \cite{LiCh} 
as well as to the appendix of \cite{GDMS24}. 
%We denote by $\delta_p$ the $p$-nuclear tensor norm, see, e.g., \cite{GDMS24}, Definition 5, and 
We use the conventions $\HH_{\alpha,1,p} := \HH_{\alpha,p}$ and $S^{\alpha}_{p}H(\real):= H^{\alpha}_{p}(\real)$.
%(Definition~\ref{tensdefi} \framebox{STILL TO DO!}-- sollten wir das auch hier im Dokument stehen haben?)\\

Important for us is the $p$-nuclear or left $p$-Chevet-Saphar norm, which is a so-called uniform cross norm.

\begin{definition}\label{tensdefi}   Let $X$ and $Y$ be Banach spaces.
Moreover, let $1 < p <  \infty$ and let $1/p+1/p' = 1$. On the algebraic tensor product $X\otimes Y$ of $X$ and $Y$ the \emph{$p$-nuclear tensor norm} 
or  \emph{left $p$-Chevet-Saphar norm} $\delta_p(\cdot,X,Y)$ is given by
\begin{equation*}
\begin{split}
    &\delta_p(h,X,Y) := \\
    &~~~\inf\left\{\Big(\sum\limits_{i=1}^n \|f_i \|^p_{X}\Big)^{1/p}\cdot
    \sup\Big\{\Big(\sum\limits_{i=1}^n |\psi(g_i)|^{p'}\Big)^{1/p'}:
    \psi \in Y', \|\psi \|_{Y'} \leq 1 \Big\}\right\}\,,
  \end{split}
\end{equation*}
where the infimum is taken over all representations of $h$ of the form
\[
h = \sum_{i=1}^n f_i \otimes g_i\, , \qquad f_i \in X\, , \quad g_i \in Y\, , \quad i=1, \ldots \, ,n\, ,
\]
and $Y'$ denotes the dual space of $Y$, i.e., the space of all continuous linear functionals on $Y$. \\
We write $X\otimes_{\delta_p}Y$ for the completion of the algebraic tensor product $X\otimes Y$ with regard to the $p$-nuclear tensor norm $\delta_p$.
\end{definition}

An alternative representation of the $p$-nuclear norm can be found in \cite[Lemma~1.44]{LiCh} or in \cite[Rem.~17]{GDMS24}.

\begin{theorem}[\cite{GDMS24}, Thm.~4]\label{tensor3}
 Let $s\ge 2$ and $1 < p<\infty$ and let $\alpha \in (1/p,1]$. 
 Then we have
 \[
 \HH_{\alpha,p}  \otimes_{\delta_p}\HH_{\alpha,s-1, p} = \HH_{\alpha,s-1,p}  \otimes_{\delta_p}\HH_{\alpha,p} =  \HH_{\alpha,s,p} 
 \]
 in the sense of equivalent norms.
\end{theorem}

A proof of Theorem~\ref{tensor3} can be found %in the appendix 
in \cite[Sect.~5.3]{GDMS24}.
Next we state a corresponding tensorization result for Bessel potential spaces.

\begin{proposition}\label{tensor4} 
Let $s\ge 2$ and let   $\alpha \in [0,1]$ and $1 < p<\infty$.
Then we have
 \begin{equation}\label{eq-5b}
H^\alpha_p([0,1]) \otimes_{\delta_p} S^\alpha_pH([0,1]^{s-1}) = S^\alpha_pH([0,1]^{s-1}) \otimes_{\delta_p} H^\alpha_p([0,1])=  S^\alpha_pH([0,1]^{s})
 \end{equation}
in the sense of equivalent norms.
\end{proposition}

A proof of Proposition \ref{tensor4} can be found in the appendix in \cite{SU11}; cf. also \cite[Prop.~5(i)]{GDMS24}.

Corollary~\ref{equal01}, Theorem \ref{tensor3}, and Proposition~\ref{tensor4}, together with a standard argument for tensor operators on tensor spaces, immediately 
result in the next corollary.

\begin{corollary}\label{equal02}
Let $1<\alpha^{-1}<p<\infty$. Then
$$S^\alpha_p H([0,1]^s)= \mathcal H_{\alpha ,s,p }$$
in the sense of equivalent norms.
\end{corollary}
This answers Question (Q2) from Section \ref{GandM}.

\subsection{Fractional Discrepancy and Koksma-Hlawka-Type Inequalities}

For a pointset $P=\{\bsx_0,\ldots,\bsx_{N-1}\}$ in $[0,1]^s$ with $\bsx_n = (x_{n,1},\ldots,x_{n,s})$
and $\emptyset\neq u\subseteq S$ we define the \emph{(local) discrepancy function $\Delta$} for $\bst_u\in [0,1]^{u}$ as
$$\Delta(\bst_u,u,P):=\prod_{j\in u}t_j-\frac 1N\sum_{n=0}^{N-1}\prod_{j\in u}\mathbb 1_{( x_{n,j},1]}(t_j);$$
in case $u = S$, we usually suppress the references to $u$ and obtain 
the classical, well-studied discrepancy function $\Delta(\cdot, P)$, anchored in $\bszero$. The \emph{($L^\infty$-)star discrepancy
of $P$} is defined as 
\[
D^*_{N}(P) := \sup_{\bsx\in [0,1)^s}|\Delta(\bsx,P)|.
\]
Let us use the short hands
\begin{equation*}
I_s(f):= \int_{[0,1]^s} f(\bst)\,{\rm d}\bst
\hspace{3ex}\text{and}\hspace{3ex}
Q_P(f) := \frac{1}{N}
\sum_{n=0}^{N-1} f(\bsx_n).
\end{equation*}
Then the \emph{classical Koksma-Hlawka inequality}, see
\cite{ Kok42, Hla61}, is of the form
\begin{equation*}
\left|I_s(f) - Q_P(f) \right|  \le  D^\ast_{N}(P)
V_{{\rm HK}}(f),
\end{equation*}
where $V_{{\rm HK}}(f)$ is the \emph{variation of $f$ in the sense of Hardy and Krause}.
%mmm

We want to study extensions of this classical error estimate for quasi-Monte Carlo methods.
Of course, many extensions to different settings are already known, see, e.g., \cite{Hic98, HSW04, NW10, BCGT2011, AD15, PS2015}.
We are mainly interested in Koksma-Hlawka inequalities that are suitable for functions with moderate fractional smoothness.
To this purpose we introduce some further notions.
For a real number $0<\alpha\le 1$, a pointset $P=\{\bsx_0,\ldots,
\bsx_{N-1}\}$ in $[0,1]^s$, 
%with $\bsx_n = (x_{n,1},\ldots,x_{n,s})$,
a set $\emptyset\neq u \subseteq S$, and a point
$\bst_u = (t_j)_{j\in u}\in [0,1]^{u}$ we define the
\emph{fractional discrepancy function} $\Delta_\alpha$ 
as in \cite{Dic08a, GDMS24}
by
\begin{equation*}
\Delta_\alpha(\bst_u, u, P) := \alpha^{-|u|}
\prod_{j\in u}(1-t_j)^\alpha - \frac{1}{N} \sum^{N-1}_{n=0}
\prod_{j\in u}(x_{n,j}-t_j)_+^{\alpha-1}.
\end{equation*}
Obviously, $\Delta_\alpha(\cdot, u, P)$ is a measurable function.
If $\alpha=1$, then the term $\prod_{j\in u}(x_{n,j}-t_j)_+^{\alpha-1}$ is, by convention, equal to $\mathbb 1_{[0,\bsx_{n,u})}(\bst_u)$, and we obtain
$\Delta_1(\cdot, u, P)\in L_\infty ([0,1]^{u})$. 
%as well as $\|\Delta_1(\cdot,u,P)\|_{L_q} >0$ for all $q\in [1,\infty]$.
If $\alpha\in (0,1)$ and $q<(1-\alpha)^{-1}$, then we have $\Delta_\alpha(\cdot, u, P)\in L_q([0,1]^{u})$.
%and $\|\Delta_\alpha(\cdot,u,P)\|_{L_q} >0$. \jo{Normen $>0$ ist interessant/wichtig? }
One can quickly see that 
$$\Delta(\bst_u,u,P)=\Delta_1(\bsone-\bst_u,u,\overline{P}),$$
 where $\overline{P}$ is the reflected pointset, i.e., for $\bsx_n\in P$ we put $\overline{\bsx}_n:=\bsone-\bsx_n\in \overline{P}$.

The following fractional Koksma-Hlawka-type inequality was shown in \cite{Dic08a} for $\alpha \in (0,1)$;
in the classical case $\alpha =1$ the result was known before, see, e.g., \cite{Hic98}.

\begin{theorem}[\cite{Dic08a}, Thm.~1]\label{th_fracKH}
Let $s \ge 1$, $0 < \alpha \le 1$, $p > \alpha^{-1}$, and $q \ge 1$.
Then for any $f \in
\HH_{\alpha,s,p}$ we have
\begin{equation*}%\label{fracKH}
\left|I_s(f) - Q_P(f) \right|  \le  D^\ast_{\alpha,s,p',q'}(P)
V_{\alpha,s,p,q}(f),
\end{equation*}
where
\begin{equation*}
D^\ast_{\alpha,s,p',q'}(P) = \left(\sum_{\emptyset \neq u \subseteq S}
 \left(\int_{[0,1]^{u}}
|\Delta_\alpha(\bst_u,u,P)|^{p'} \,\rd\bst_u \right)^{q'/p'}
\right)^{1/q'}
\end{equation*}
and
\begin{equation*}
V_{\alpha,s,p,q}(f) =\left(\sum_{\emptyset \neq u \subseteq S} 
\Gamma(\alpha)^{-q|u|}
\left(\int_{[0,1]^{u}} \left| D^{(\bsalpha_u, \bszero_{S\setminus u})}f(\bst_u, \bszero)\right|^{p} \,\rd\bst_u
\right)^{q/p} \right)^{1/q}
\end{equation*}
with the obvious modifications in the case where $p,q$ or $q'=\infty$.
\end{theorem}

\begin{remark}\label{Rem:Seminorm_Norm}
Note that
$V_{\alpha,s,p,q}(\cdot)$ is just the seminorm in the space
$\HH_{\alpha,s,p}$ and could be replaced by the norm
$\|\cdot\|_{\alpha,s,p,q}$.
Furthermore, for $\alpha\in (0,1)$ the condition $p > \alpha^{-1}$ implies $p' < (1-\alpha)^{-1}$, which
is also the condition required to guarantee that
$D^\ast_{\alpha,s,p',q'}(P) < \infty$. 
We refer to $D^\ast_{\alpha,s,p',q'}(P)$ as \emph{fractional $(p', q')$-discrepancy} 
of the pointset $P$.
\end{remark}

The fractional Koksma-Hlawka inequality is indeed
sharp. To formulate this fact in the form of an identity, let us denote by $e^{\wor}(Q_P,\HH_{\alpha,s,p,q})$ the
worst-case error
of the quasi-Monte Carlo quadrature $Q_P$ based on the set of
integration points $P$ on $\HH_{\alpha,s,p,q}$,  i.e.,
\begin{equation*}
e^{\wor}(Q_P,\HH_{\alpha,s,p,q}) = \sup_{f\in \HH_{\alpha,s,p}; \|f\|_{\alpha,s,p,q}\le 1} \left| \int_{[0,1]^s} f(\bsx) \, \rd \bsx
- \frac{1}{N} \sum_{n=0}^{N-1} f(\bsx_n) \right|.
\end{equation*}
The worst-case error $e^{\wor}(Q_P,S^{\alpha}_{p}H([0,1]^s))$
of $Q_P$ on the Bessel potential space $S^{\alpha}_{p}H([0,1]^s)$ is defined analogously.
The following theorem holds.

\begin{theorem}[\cite{GDMS24}, Thm.~6]\label{KHI-GDMS}
\label{KoksmaHlawka}
Let $s\in\NN$, $0<\alpha\le 1$, $p>\alpha^{-1}$ and $q\in [1,\infty]$.
For all pointsets $P\subseteq [0,1]^s$ we have
\begin{equation}
\label{koksmahlawka}
e^{\wor}(Q_P,\HH_{\alpha,s,p,q}) =  D^\ast_{\alpha,s,p',q'}(P).
\end{equation}
\end{theorem}

In the case $\alpha = 1$ this result was known earlier, cf., e.g., \cite[Sect.~9.8.1]{NW10}.

Theorem~\ref{KoksmaHlawka} leads to the following insight: 
If we rely on the fractional discrepancy as quality criterion for sample points, then the spaces of fractional smoothness are the appropriate function spaces for the analysis of QMC quadrature rules.\\

\noindent
We now use these fractional Koksma-Hlawka inequalites to obtain Koksma-Hlawka type inequalites for Bessel potential spaces. 
We begin by stating a direct consequence of Corollary \ref{equal02} and Theorem \ref{KHI-GDMS}.

\begin{corollary}[First Koksma-Hlawka-Type Inequality for Bessel Potential Spaces]\label{KHI1}
Let $\alpha \in (0,1)$, $\alpha^{-1}<p<\infty$, and $q\ge 1$. Then there exist constants $C,c>0$ such that for all pointsets $P=\{\bsx_0,...,\bsx_{N-1}\}\subseteq[0,1]^s$ we have
\begin{equation*}
c D^\ast_{\alpha,s,p',q}(P) \le
%\sup_{ ||f||_{S^{\alpha}_{p}H([0,1]^s)}\leq1}\left|I_s(f) - Q_P(f) \right|  
e^{\wor}(Q_P,S^{\alpha}_{p}H([0,1]^s))
\le C D^\ast_{\alpha,s,p',q}(P).
\end{equation*}
\end{corollary}

In \cite[Sect.~6]{Tri10} Koksma-Hlawka type inequalities have been proved for Besov spaces with boundary conditions. Instead of applying an $L^p$-norm, these inequalities use the norm of the discrepancy function in an adjoint Besov space. We now want to obtain similar results for Bessel potential spaces, 
%we however don't need the boundary conditions. 
but preferably without having to assume any addtional boundary conditions. We start with a crucial observation.

\begin{lemma}\label{Lemma:CrucialObs}
Let $s\in\N$ and $1<\alpha^{-1}<p<\infty$. Let furthermore $\emptyset\neq u\subseteq S$ and $P=\{\bsx_0,...,\bsx_{N-1}\}$ be a pointset in $[0,1]^s$.
Define $\tilde\Delta:\R^u\to\R$ by
$$\tilde \Delta:=\Delta^0(\cdot,u,\overline{P})*\prod_{j\in u}(\cdot_j)^{\alpha-1}_+,$$
where $\Delta^0(\cdot,u,\overline P)$ is the $0$-extension of $\Delta(\cdot,u,\overline P)$ on $\R^u$.
\begin{itemize}
    \item [(i)] For $|u|=n\in\N$, $u=\{s_1,...,s_n\}$, we use the notation 
$\frac{\partial^{|u|}}{\partial \bst^u}:=\frac{\partial}{\partial t_{s_1}}...\frac{\partial}{\partial t_{s_n}}.$
    It is
$$\frac{\partial^{|u|}}{\partial \bst^u}\tilde \Delta(\bst_u)=\Delta_\alpha(\bsone-\bst_u,u,P)$$
%\big(\Delta^0(\cdot,u,\overline{P})*\prod_{j\in u}(\cdot_j)^{\alpha-1}_+\big)
for almost all $\bst_u\in[0,1]^u$. 
\item[(ii)] We have $\tilde \Delta|_{[-1,2]^u} \in S^1_{p'}H([-1,2]^u)$ and its mixed derivative almost everywhere from $(i)$ coincides for almost all $\bst_u\in[0,1]^u$ with the corresponding distributional derivative on $[-1,2]^u$. 
\end{itemize}
\begin{comment}
Let $s\in\N$ and let $1<\alpha^{-1}<p<\infty$. Then for all $\emptyset\neq u\subseteq S$ and for each pointset $P=\{\bsx_0,...,\bsx_{N-1}\}$ in $[0,1]^s$ we have
$$\Delta_\alpha(\bsone-\bst_u,u,P)=\frac{\partial^{|u|}}{\partial \bst^u}\big(\Delta^0(\cdot,u,\overline{P})*\prod_{j\in u}(\cdot_j)^{\alpha-1}_+\big)(\bst_u)$$
%=\frac{\partial^{|u|}}{\partial \bst^u}\int_{[0,1]^u}\Delta(\bsx_u,u,\bar{P})\prod_{j\in u}(\bst_j-\bsx_j)^{\alpha-1}_+\,d\bsx_u$$ 
for almost all $\bst_u\in[0,1]^u$, where $\Delta^0(\cdot,u,\overline P)$ is the $0$-extension of $\Delta(\cdot,u,\overline P)$ on $\R^u$.    
\end{comment}
\end{lemma}

\begin{proof}
 % We show 
%$\Delta_\alpha(\bsone-\bst_u,u,P)=\frac{\partial^{|u|}}{\partial t^u}\big(\Delta^0(\cdot,u,\overline{P})*\prod_{j\in u}(\cdot_j)^{\alpha-1}_+\big)(\bst_u)$ for almost all $\bst_u\in[0,1]^u$, where $\Delta^0$ is the $0$-continuation of $\Delta$ on $\R^u$.\\
We begin with a calculation needed for the proofs of both parts of the Lemma.\\
Recall that
%$$\Delta(t_u,u,P):=\prod_{j\in u}t_j-\frac 1N\sum_{n=0}^{N-1}\prod_{j\in u}\mathbb 1_{( x_{n,j},1]}(t_j)$$
% and 
\begin{equation*}
\Delta_\alpha(\bst_u, u, P) := \alpha^{-|u|}
\prod_{j\in u}(1-t_j)^\alpha - \frac{1}{N} \sum^{N-1}_{n=0}
\prod_{j\in u}(x_{n,j}-t_j)_+^{\alpha-1}.
\end{equation*}
In the following we use the notation $x\wedge y$ and $x\vee y$ to denote the minimum and the maximum of two real numbers $x$ and $y$.
For $\bst_u\in\R^u$ we get 
%\framebox{NOTATION $\vee$ UND $\wedge$ NOCH EINFUEHREN!}
     \begin{eqnarray}\label{tensorstr}
    \begin{aligned}
    &\big(\Delta^0(\cdot,u,\overline{P})*\prod_{j\in u}(\cdot_j)^{\alpha-1}_+\big)(\bst_u)\\
    =&\int_{[0,1]^u}\Big(\prod_{j\in u}y_j-\frac1N\sum_{n=0}^{N-1}\prod_{j\in u}\mathbb 1_{( \overline x_{n,j},1]}(y_j)\Big)\prod_{j\in u}(t_j-y_j)^{\alpha-1}_+\,d\bsy_u\\
    =&\prod_{j\in u}\int_0^{0 \vee (1 \wedge t_j)}y_j(t_j-y_j)^{\alpha-1} \,dy_j-\frac1N\sum_{n=0}^{N-1}\prod_{j\in u}
    \int_{\overline x_{n,j}}^{\overline x_{n,j} \vee (1 \wedge t_j)}(t_j-y_j)^{\alpha-1} \,dy_j\\
    =&\prod_{j\in u}\int_0^{0 \vee (1 \wedge t_j)} y_j(t_j-y_j)^{\alpha-1} \,dy_j-\frac1N\sum_{n=0}^{N-1}\prod_{j\in u}\frac{-1}{\alpha}\Big(\big(t_j-
    (\overline x_{n,j} \vee (1 \wedge t_j))\big)^{\alpha} - (t_j- \overline x_{n,j})^\alpha\Big)\\
    =&\prod_{j\in u} \frac{(t_j-y_j)^\alpha (t_j+\alpha y_j)}{-\alpha(\alpha+1)}\Big|_0^{0 \vee (1 \wedge t_j)}-\frac1N\sum_{n=0}^{N-1}\Big(\frac{-1}{\alpha}\Big)^{|u|}\prod_{j\in u}\Big(\big(t_j- (\overline x_{n,j} \vee (1 \wedge t_j))\big)^{\alpha}-(t_j- \overline x_{n,j})^\alpha\Big)\\
    =&\prod_{j\in u} \frac{(t_j)^{\alpha+1}- \big(t_j- (0 \vee (1 \wedge t_j))\big)^\alpha \big(t_j+\alpha (0 \vee (1 \wedge t_j)) \big)}{\alpha(\alpha+1)}\\
    -&\frac1N\sum_{n=0}^{N-1}\Big(\frac{-1}{\alpha}\Big)^{|u|}\prod_{j\in u}\Big(\big(t_j- (\overline x_{n,j} \vee (1 \wedge t_j))\big)^{\alpha}-(t_j- \overline x_{n,j})^\alpha\Big),\\
    \end{aligned}
    \end{eqnarray}
    where we use the convention, that for $x<0$ and $p>0$  one has $x^p-x^p=0$, even though these powers aren't formally defined.\\
    $(i)$: We see that $\big(\Delta^0(\cdot,u,\overline{P})*\prod_{j\in u}(\cdot_j)^{\alpha-1}_+\big)(\bst_u)$ is continuous in $\bst_u$ 
   % and that its restriction to $[-1,2]^u$ is a finite sum of tensor products of functions in $H^1_{p'}([-1,2])$, hence it lies itself in $S^1_{p'}H([-1,2]^u)$. \\
    %(First summand: $t_j$ less than $1$: Beta function and cont. $t_j>1$ cont. since there is not even a singularity anymore.)\\
    and get for almost all $\bst_u\in [0,1]^u$
      \begin{eqnarray*}
    \begin{aligned}
    &\frac{\partial^{|u|}}{\partial \bst^u}\big(\Delta^0(\cdot,u,\overline{P})*\prod_{j\in u}(\cdot_j)^{\alpha-1}_+\big)(\bst_u)\\
    =&\frac{\partial^{|u|}}{\partial \bst^u}\Big(\prod_{j\in u}\frac{t_j^{\alpha+1}}{\alpha(\alpha+1)}-\frac1N\sum_{n=0}^{N-1}\frac{1}{\alpha^{|u|}}\prod_{j\in u} (t_j- \overline x_{n,j})^\alpha_+\Big)\\
    =&\frac1{\alpha^{|u|}}\prod_{j\in u}t_j^{\alpha}-\frac1N\sum_{n=0}^{N-1}\prod_{j\in u} (t_j- \overline x_{n,j})^{\alpha-1}_+\\
    =&\Delta_\alpha(\bsone-\bst_u,u,P).
    \end{aligned}
    \end{eqnarray*}
  %  We obviously also have the equality for both the $0$-continuations on $\R^u$.
  \begin{comment}
       \jo{
    Furthermore, this derivative almost everywhere coincides with the distributional derivative of $\Delta^0(\cdot,u,\overline{P})*\prod_{j\in u}(\cdot_j)^{\alpha-1}_+$ on $[0,1]^u$.
    The reason for this is that $\Delta^0(\cdot,u,\overline{P})*\prod_{j\in u}(\cdot_j)^{\alpha-1}_+$ is a finite sum of tensor products of absolutely continuous functions on 
    $[-1,2]^u$ and therefore absolutely continuous along all lines parallel to a coordinate axis. It is known, that for absolutely continuous functions on intervals the distributional derivative and the derivative almost everywhere coincide, see e.g. \cite{Rud73}. We can therefore use this statement iteratively in every dimension to get the equality of the different derivatives for $\Delta^0(\cdot,u,\overline{P})*\prod_{j\in u}(\cdot_j)^{\alpha-1}_+$.
    }
  \end{comment}
  $(ii)$: We proceed similarly to \cite[Thm. 2.1.4]{Zi89}. Using the tensor structure of $\tilde \Delta$ established in \eqref{tensorstr}, one can quickly show that $\tilde \Delta$ is absolutely continuous along almost all line segments parallel to the coordinate axes in $[-1,2]^u$. 
For every $\phi\in C^\infty_0([-1,2]^u)$ the same property holds for $\phi\cdot \tilde\Delta$.
%Thus, since for absolutely continuous functions on intervals the distributional derivative and the derivative almost everywhere coincide, see e.g. \cite{Rud73}, it follows that for all $i\in u$ one has 
Thus, by Fubini's theorem and Lebesgue's version of the fundamental theorem of calculus, one obtains for each $i\in u$  
$$\int_{[-1,2]^u}\phi(\bst_u)\frac{\partial}{\partial t_i}\tilde \Delta(\bst_u) \, d\bst_u=-
\int_{[-1,2]^u}\big(\frac{\partial}{\partial t_i}\phi\big)(\bst_u)\tilde \Delta(\bst_u) \, d\bst_u,$$
which implies
$$\frac{\partial}{\partial t_i}\tilde \Delta(\bst_u)=
D_i\tilde \Delta(\bst_u)$$
for almost all $\bst_u\in[-1,2]^u$, where by $D_i$ we denote the $i$-th distributional derivative in $[-1,2]^u$. 
Now for $j\in u$ with 
$j\neq i$ we have that $\frac{\partial}{\partial t_i}\tilde \Delta$ is absolutely continuous along almost all line segments in $[-1,2]^u$ parallel to the $j$-axis, implying as before
$$\frac{\partial^2}{\partial t_j\partial t_i}\tilde \Delta(\bst_u)=
D_jD_i\tilde \Delta(\bst_u)$$
for almost all $\bst_u\in[-1,2]^u$. 
The proof is finished by using this procedure iteratively and the fact that all these derivatives lie in $L^{p'}([-1,2]^u)$.
%and Theorem \ref{BessSob} 
%we obtain $\tilde \Delta \in S^1_{p'}H([-1,2]^u)$ and that its mixed derivative almost everywhere from Lemma \ref{Lemma:CrucialObs}
% coincides with the corresponding distributional derivative for almost all $\bst_u\in[0,1]^u$. 
\end{proof}

\begin{comment}
\jo{
\begin{remark}
For $s=1$, $1<\alpha^{-1}<p<\infty$ and pointsets $P=\{x_0,...,x_{N-1}\}\subseteq[0,1]$, the statement of Lemma \ref{Lemma:CrucialObs} becomes 
$$\Delta_\alpha(1-t,P)= \Gamma(\alpha)\frac{d^{1-\alpha}}{dt^{1-\alpha}}\Delta(t,\overline{P}).$$ 
With this, the $L^p$-norm of the fractional discrepancy function $\Delta_\alpha$ looks similar to the norm of the discrepancy function $\Delta$ in a space of fractional smoothness. We even have due to \cite[Prop. 6.3]{Tri10} 
    $$\Delta(\cdot,\overline{P})\in H^{1-\alpha}_{p'}([0,1]),$$
since $(1-\alpha)^{-1}>p'$. With this and Corollary \ref{KHI1} it would be natural to use Corollary \ref{equal02} in order to get a Koksma-Hlawka type inequality using norms of Bessel potential spaces on the discrepancy function and the integrated function. The problem however is, that for the parameter range $(1-\alpha)^{-1}>p'$ the spaces of fractional smoothness $\mathcal H_{1-\alpha,p'}$ are not defined and therefore Corollary \ref{equal02} can not be applied.
\end{remark}
% In order to obtain Koksma-Hlawka inequalities using the norms of two adjoint Bessel potential spaces on the integrated function and the discrepancy function, we now use methods from harmonic analysis.
 In order to still obtain this type of Koksma-Hlawka inequality, we will use similar methods from harmonic analysis as in the proof of Corollary \ref{equal01}. We 
 will however work explicitly with the discrepancy function and not an arbitrary function from a function space.
}    
\end{comment}
 In order to obtain Koksma-Hlawka inequalities using the norms of two adjoint Bessel potential spaces on the integrated function and the discrepancy function, we now use similar methods from harmonic analysis as in the proof of Theorem \ref{embedding01}. We will however work explicitly with the discrepancy function and not an arbitrary function from a function space.

\begin{theorem}\label{BessDisc}
    Let $s\in\N$ and let $1<\alpha^{-1}<p<\infty$. Then for all $\emptyset\neq u\subseteq S$ there exists a constant $C>0$ such that for each pointset $P=\{\bsx_0,...,\bsx_{N-1}\}$ in $[0,1]^s$ we have
    $$||\Delta_\alpha(\cdot,u,P)||_{L^{p'}([0,1]^u)}
    \leq C||\Delta_1(\cdot,u,P)||_{S^{1-\alpha}_{p'}H([0,1]^u)}.$$ 
\end{theorem}

\begin{proof}
% The theorem is proved in two steps.\\
%\textbf{Step 1:} 
%\textbf{Step 2:} We finish the proof by using tools from harmonic analysis.\\
    Due to \cite[Prop. 6.3]{Tri10} we get
    $$\Delta(\cdot,u,\overline{P})\in S^{1-\alpha}_{p'}H([0,1]^u),$$
    since $(1-\alpha)^{-1}>p'$. This implies that there exists an extension $\bar \Delta(\cdot,u,\overline{P})\in S^{1-\alpha}_{p'}H(\R^u)$ of $\Delta(\cdot,u,\overline{P})$ with 
    $$||\bar\Delta(\cdot,u,\overline{P})||_{S^{1-\alpha}_{p'}H(\R^u)}\leq 2||\Delta(\cdot,u,\overline{P})||_{S^{1-\alpha}_{p'}H([0,1]^u)}.$$
    Since due to \cite[Cor. II3.7]{Str67} the characteristic function of a univariate interval is a pointwise multiplier on $H^{1-\alpha}_{p'}(\R)$ and due to the tensor structure of $S^{1-\alpha}_{p'}H(\R^u)$,
    we get that characteristic functions of axes-parallel rectangles are pointwise multipliers on $S^{1-\alpha}_{p'}H(\R^u)$. Therefore, by applying the characteristic function of the unit cube as a pointwise multiplier to $\bar\Delta$, it follows that there is a constant $\tilde C>0$ such that for the 
    $0$-extension $\Delta^0(\cdot,u,\overline{P})$ of $\Delta(\cdot,u,\overline{P})$ onto $\R^u$ we have
    $$\tilde C||\Delta^0(\cdot,u,\overline{P})||_{S^{1-\alpha}_{p'}H(\R^u)}\leq ||\Delta(\cdot,u,\overline{P})||_{S^{1-\alpha}_{p'}H([0,1]^u)}.$$
    With the function $m$ from Example \ref{ex2} we now get 
    \begin{comment}
          Now put
    $$m_0(\xi):=\prod_{j\in u}\sign(\xi_j)e^{-i\pi\frac\alpha 2 \sign(\xi_j)}$$
    and for all $j\in u$
    $$m_j(\xi):=\frac{|\xi_j|^{1-\alpha}}{<\xi_j>^{1-\alpha}}.$$
    By Corollary \ref{MarMT} it is easily seen, that these are $L^p$-multipliers for $1<p<\infty$. Therefore this also holds for $m:=\prod_{j\in u\cup\{0\}}m_j$.
    We now get
    \end{comment}
    \begin{eqnarray*}
    \begin{aligned}
    ||\Delta^0(\cdot,u,\overline{P})||_{S^{1-\alpha}_{p'}H(\R^u)}=&||\mathcal F^{-1} \prod_{j\in u}<\xi_j>^{1-\alpha}\mathcal F\big(\Delta^0(\cdot,u,\overline{P})\big)||_{L^{p'}(\R^u)}\\
    \geq&C ||\mathcal F^{-1}\Big(m(\xi) \prod_{j\in u}<\xi_j>^{1-\alpha}\mathcal F\big(\Delta^0(\cdot,u,\overline{P})\big) \Big)||_{L^{p'}(\R^u)}.\\
    \end{aligned}
    \end{eqnarray*}
    
    Similarly to the univariate case, we introduce mollifiers. For a functions $\psi\in C^\infty_0(\R^s)$ with $\int_{\R^s} \psi(\bsx)\,d\bsx=1$, $\psi(\bsx)\geq0$ for all $\bsx\in\R^s$ and $\supp(\psi)\subset [-1,1]^s$, we call the sequence
    of functions $(\psi_n)_n$ with $\psi_n(\bsx):=n^s\psi(n\bsx)$ a mollifier, see e.g. \cite{Gra14}.\\
    Let now $(\psi_n^u)_n\subseteq C^\infty_0(\R^u)$ be a mollifier. \\
    By using the convolution theorem for tempered distributions twice (see e.g. \cite[Prop. 2.3.22]{Gra14}), we get
     \begin{eqnarray*}
    \begin{aligned}
    &||\mathcal F^{-1}\Big(m(\xi) \prod_{j\in u}<\xi_j>^{1-\alpha}\mathcal F\big(\psi_n^u*\Delta^0(\cdot,u,\overline{P})\big) \Big)||_{L^{p'}(\R^u)}\\
    =&||\psi_n^u*\mathcal F^{-1}\Big(m(\xi)\prod_{j\in u}<\xi_j>^{1-\alpha}\mathcal F\big(\Delta^0(\cdot,u,\overline{P})\big) \Big)||_{L^{p'}(\R^u)}\\
    \to&||\mathcal F^{-1}\Big(m(\xi) \prod_{j\in u}<\xi_j>^{1-\alpha}\mathcal F\big(\Delta^0(\cdot,u,\overline{P})\big) \Big)||_{L^{p'}(\R^u)}
    \end{aligned}
    \end{eqnarray*}
    if $n$ tends to infinity.
    However, we also get
    \begin{eqnarray*}
    \begin{aligned}
    &||\mathcal F^{-1}\Big(m(\xi) \prod_{j\in u}<\xi_j>^{1-\alpha}\mathcal F\big(\psi_n^u*\Delta^0(\cdot,u,\overline{P})\big) \Big)||_{L^{p'}(\R^u)}\\
    =&||\mathcal F^{-1}\Big(\prod_{j\in u}\frac{|\xi_j|^{1-\alpha}}{<\xi_j>^{1-\alpha}}\sign(\xi_j)e^{-i\pi\frac\alpha 2 \sign(\xi_j)} \prod_{j\in u}<\xi_j>^{1-\alpha} \mathcal F\big(\psi_n^u*\Delta^0(\cdot,u,\overline{P})\big) \Big)||_{L^{p'}(\R^u)}\\
     =&||\mathcal F^{-1}\Big(\prod_{j\in u}|\xi_j|^{1-\alpha}\sign(\xi_j)e^{-i\pi\frac\alpha 2 \sign(\xi_j)}\cdot \mathcal F\big(\psi_n^u*\Delta^0(\cdot,u,\overline{P})\big)\Big)||_{L^{p'}(\R^u)}\\
    =&||\mathcal F^{-1}\Big(\prod_{j\in u}\xi_j\cdot\mathcal F\big(\psi_n^u*\Delta^0(\cdot,u,\overline{P})\big)\cdot \prod_{j\in u}%\frac{1}{\sqrt(2\pi)^{|u|}}
    |\xi_j|^{-\alpha}e^{-i\pi\frac\alpha 2 \sign(\xi_j)}\Big)||_{L^{p'}(\R^u)}\\
    =&c||\mathcal F^{-1}\Big(\mathcal F\frac{\partial^{|u|}}{\partial t^u}\big(\psi_n^u*\Delta^0(\cdot,u,\overline{P})\big)\cdot\mathcal F\prod_{j\in u}(\cdot_j)^{\alpha-1}_+\Big)||_{L^{p'}(\R^u)},
    \end{aligned}
    \end{eqnarray*}
    where in the last step we used the well-known interactions between the Fourier-transform and differentiation, as well as 
    $$\F((\cdot)_+^{\alpha-1})(\xi)=\frac{\Gamma(\alpha)}{\sqrt{2\pi}}|\xi|^{-\alpha}e^{-i\pi \frac{\alpha}{2}\sgn(\xi)}$$
    from \cite[p. 400]{BN71}; the constant factors arising in the process form the constant $c=\Gamma(\alpha)^{|u|}$.\\
    The convolution theorem for distributions and some well-known results on distributional derivatives of convolutions (see e.g. \cite{Fo99, Gra14}) now imply
    \begin{eqnarray*}
    \begin{aligned}
    &c||\mathcal F^{-1}\Big(\mathcal F\frac{\partial^{|u|}}{\partial t^u}\big(\psi_n^u*\Delta^0(\cdot,u,\overline{P})\big)\cdot\mathcal F\prod_{j\in u}(\cdot_j)^{\alpha-1}_+\Big)||_{L^{p'}(\R^u)}\\
    =&c||\mathcal F^{-1}\mathcal F\Big(\frac{\partial^{|u|}}{\partial t^u}\big(\psi_n^u*\Delta^0(\cdot,u,\overline{P})\big)*\prod_{j\in u}(\cdot_j)^{\alpha-1}_+\Big)||_{L^{p'}(\R^u)}\\
    =&c||\frac{\partial^{|u|}}{\partial t^u}\big(\psi_n^u*\Delta^0(\cdot,u,\overline{P})\big)*\prod_{j\in u}(\cdot_j)^{\alpha-1}_+||_{L^{p'}(\R^u)}\\
    =&c||\big(\frac{\partial^{|u|}}{\partial t^u}\psi_n^u*\Delta^0(\cdot,u,\overline{P})\big)*\prod_{j\in u}(\cdot_j)^{\alpha-1}_+||_{L^{p'}(\R^u)}\\
    =&c||\frac{\partial^{|u|}}{\partial t^u}\psi_n^u*\big(\Delta^0(\cdot,u,\overline{P})*\prod_{j\in u}(\cdot_j)^{\alpha-1}_+\big)||_{L^{p'}(\R^u)}\\
    \geq& c||\frac{\partial^{|u|}}{\partial t^u}\psi_n^u*\big(\Delta^0(\cdot,u,\overline{P})*\prod_{j\in u}(\cdot_j)^{\alpha-1}_+\big)||_{L^{p'}([0,1]^u)}\\
    =&c||\psi_n^u*\frac{\partial^{|u|}}{\partial t^u}\big(\Delta^0(\cdot,u,\overline{P})*\prod_{j\in u}(\cdot_j)^{\alpha-1}_+\big)||_{L^{p'}([0,1]^u)}\\
    \to&c||\frac{\partial^{|u|}}{\partial t^u}\big(\Delta^0(\cdot,u,\overline{P})*\prod_{j\in u}(\cdot_j)^{\alpha-1}_+\big)||_{L^{p'}([0,1]^u)}
    %\overset{**}{=}&c||\Delta_\alpha(\bsone-\cdot,u,P)||_{L^{p'}([0,1]^u)}.
    %=||\Delta_\alpha(\cdot,u,P)||_{L^{p'}([0,1]^u)}.
    %\geq& c||\psi_n^u*\frac{\partial^{|u|}}{\partial t^u}\big(\Delta^0(\cdot,u,P)*\prod_{j\in u}(\cdot_j)^{\alpha-1}_+\big)||_{L^{p'}([0,1]^u)}
    %\to& c||\frac{\partial^{|u|}}{\partial t^u}\big(\Delta^0(\cdot,u,P)*\prod_{j\in u}(\cdot_j)^{\alpha-1}_+\big)||_{L^{p'}(\R^u)}.
    \end{aligned}
    \end{eqnarray*}
    if $n$ tends to infinity. 
    Note that in the penultimate step above it is eligible to shuffle the partial derivatives from the first factor
    of the convolution into the second one. This is due to the second factor being in $S^1_{p'}H([-1,2]^u)$, 
    see Lemma \ref{Lemma:CrucialObs} $(ii)$, and   \cite[Prop. 9.3]{Fo99}. 
    Observe that Lemma~\ref{Lemma:CrucialObs} implies that
    \begin{equation*}
    ||\frac{\partial^{|u|}}{\partial t^u}\big(\Delta^0(\cdot,u,\overline{P})*\prod_{j\in u}(\cdot_j)^{\alpha-1}_+\big)||_{L^{p'}([0,1]^u)}
    = ||\Delta_\alpha(\bsone-\cdot,u,P)||_{L^{p'}([0,1]^u)}. 
    \end{equation*}
    % since the function is continuous (check directly) and a finite sum of tensor products of functions from $H^1_{p'}([0,1]^u)$ 
    %(these have mainly one dent).\\
   % This finishes the proof.
   This implies that there exists the desired constant $C>0$ with
   $$||\Delta_\alpha(\cdot,u,P)||_{L^{p'}([0,1]^u)}\leq C ||\Delta(\bsone-\cdot,u,\overline{P})||_{S^{1-\alpha}_{p'}H([0,1]^u)}= 
   C ||\Delta_1(\cdot,u,P)||_{S^{1-\alpha}_{p'}H([0,1]^u)}.$$
   \begin{comment}
     The interactions between the Fourier transform and shifts and reflections of functions are well understood. 
   Hence for $f\in S^{1-\alpha}_{p'}H(\R^u)$ one can quickly show
   $$||f(\cdot)||_{S^{1-\alpha}_{p'}H(\R^u)}=||f(\bsone-\cdot)||_{S^{1-\alpha}_{p'}H(\R^u)}.$$
   This implies the statement.  
   \end{comment}
\end{proof}

Corollary~\ref{equal02}, Theorem~\ref{th_fracKH}, Remark~\ref{Rem:Seminorm_Norm}, and Theorem~\ref{BessDisc}
now directly imply the following alternative Koksma-Hlawka-type inequality for Bessel potential spaces.

\begin{corollary}[Second Koksma-Hlawka-Type Inequality for Bessel Potential Spaces]\label{KHI2}
Let $1<\alpha^{-1}<p<\infty$ and $q\ge 1$. Then there exists a constant $C>0$ such that for each pointset $P$ in $[0,1]^s$ 
and each $f \in S^{\alpha}_{p}H([0,1]^s)$ we have
\begin{equation*}%\label{BessKH}
\left|I_s(f) - Q_P(f) \right|  \le C D^{BP}_{1-\alpha,s,p',q}(P)
||f||_{S^{\alpha}_{p}H([0,1]^s)},
\end{equation*}
where
\begin{equation*}
D^{BP}_{1-\alpha,s,p',q}(P) := \left(\sum_{\emptyset \neq u \subseteq S}
 ||\Delta_1(\cdot,u,P)||^{q}_{S^{1-\alpha}_{p'}H([0,1]^u)}
\right)^{1/q}.
\end{equation*}
\end{corollary}
This answers Question (Q3) from Section \ref{GandM}.

\subsection{Optimal Cubature Rules}
Let $b\in\N_{\geq 2}$, $s,m\in\N$ and $t\in \{0,1,...,m\}$. A $(t,m,s)$-net in base $b$ is a specific pointset, consisting of $b^m$ points, regular in the $s$-dimensional unit cube.
For a definition, see \cite{GDMS24}, Section 2.3, or the standard references \cite{Nie92} and \cite{DP10}.\\
Let us first restate the upper error bounds from \cite{GDMS24}.
\begin{corollary}[\cite{GDMS24}, Cor. 8]\label{Osc} \label{RatesGDMS}
Let $s\in\N$, $1 < p < \infty$, $1\le q \le \infty$ and  $\alpha \in (1/p, 1)$. 
Let $m\ge 1$. Let $P$ be a $(t,m,s)$-net in base $b$ and denote the corresponding QMC-cubature rule by $Q_P$.\\
%(cf.  \eqref{qmc-alg}). ???\\
Furthermore, put $N:= |P| = b^m$.
\begin{itemize}
\item[{\rm (}i{\rm )}] Let additionally $\alpha > 1/2$. Then there exists a constant $C$, depending only on $\alpha$,
$s$, $p$, $q$, $b$, and $t$, such that
\begin{equation}\label{frac_disc_estimate}
e^{\wor}(Q_P, \HH_{\alpha,s,p, q})  
\le C N^{-\alpha} \ln(N)^{\frac{s-1}{2}}.
\end{equation}
\item[{\rm (}ii{\rm )}] Let additionally $\alpha \le 1/2$. Then for each $\widetilde{p} \in (\alpha^{-1}, p]$ there exists a constant $C$, depending only on $\alpha$,
$s$, $p$, $\widetilde{p}$, $q$, $b$, and $t$, such that
\begin{equation}\label{frac_disc_estimate_2}
e^{\wor}(Q_P, \HH_{\alpha,s,p, q}) 
\le C N^{-\alpha} \ln(N)^{\frac{s-1}{{\widetilde{p}}'}}.
\end{equation}
\end{itemize}
\end{corollary}

\begin{remark}\label{Rem:Known_QMC_Int_Disc} 
QMC quadrature rules have the advantage that they are equal-weight rules and that their weights sum up to one. 
%In \cite{Mar13} it was shown that the upper error bound~\eqref{nets_besov} 
%holds for the Chen-Skriganov construction of $(t,m,s)$-nets $P$ for large base $b$, see~\cite{CS02, Skr06}, in the case where $1\le p,q \le \infty$ and $1/p < \alpha <1$.  \\
%Note that we proved that the bound~\eqref{nets_besov} holds actually for arbitrary $(t,m,s)$-nets $P$, regardless of the size of the base $b$ or of the specific construction principle.
\end{remark}

%\begin{comment}
    Using the lower error-bounds for functions in $S^\alpha_pH(\R^s)$ with support in the unit cube as discussed in \cite{UU15} and the continuity of the restriction onto $[0,1]^s$, we obtain lower error bounds in $S^\alpha_pH([0,1]^s)$. Similar lower error-bounds for periodic Bessel potential spaces have been proved in \cite{Tem15}. Combining these lower bounds with the upper error-bounds from \cite{GDMS24} for fractional spaces in the sense of Riemann-Liouville (see Corollary \ref{RatesGDMS}) and Corollary \ref{equal02}, we obtain the following result.
    
    \begin{corollary}\label{Hauptresultat}
    Let $s\in\N$, $p\in (1,\infty)$ and $\alpha>\max\{1/p,1/2\}$. 
   Then the upper bound in \eqref{frac_disc_estimate} is  best possible and for fixed $t\in\N$ %that QMC rules based on 
arbitrary $(t,m,s)$-nets $P$ in base $b$ achieve this bound. %with a constant $C$ depending only on $b,t,s$ as well as on $p$ and $\alpha$. 
\end{corollary}
This answers Question (Q1) from Section \ref{GandM}.

\section{Appendix: Function Expansions}

%\subsection{Function Extensions}
We want to discuss the extension theorem \cite[Thm. 4.5.2]{Tri92} referenced in the proof of Lemma \ref{Lemma1}. 
We give a short theoretical introduction and afterwards a version of the statement that is tailored to our needs. Since the proof in the original work is rather brief, we provide a more elaborate version for the convenience of the reader. For the theorem in its full generality, see \cite{Tri92}.

Let $\Omega= \R$ or $\Omega =\R_+$ and $B(x,t):=\{y\in\Omega:\, |x-y|< t\},\, x\in\Omega,\,t>0$. Let $|B(x,t)|$ be the length of $B(x,t)$. We write
$$\intbar_{B(x,t)}g(y)\,dy:=|B(x,t)|^{-1}\int_{B(x,t)}g(y)\,dy$$
for the mean value. Let $M\in \N\cup {\{-1\}}$, then $\mathcal P_M$ stands for the collection of all polynomials (with complex coefficients) of degree less than or equal to $M$, where we put 
$\mathcal P_{-1}=\{0\}$.
\begin{definition}[\cite{Tri92}, Def. 3.4.1]\label{Osc} 
    Let $\Omega= \R$ or $\Omega=\R_+$ and $M\in \N \cup {\{-1\}}$. Then for 
    $f\in L^1(\Omega)$ the quantity
    $$\osc^Mf(x,t):=\inf_{p\in\mathcal P_M}\intbar_{B(x,t)}|f(y)-p(y)|\,dy$$
    is called (local) oscillation.
\end{definition}
Of importance will be the following equivalent characterization of Bessel potential spaces.
\begin{theorem}[\cite{Tri92}, Thm. 3.5.1]\label{Osceq}
    Let $1\leq p<\infty,\,\alpha>0$ and $M\in\N$ with $M>\alpha$. Then
    $$H^\alpha_p(\R)=\Big\{f\in L^p(\R)\,\Big|\,||f||_{H^\alpha_p(\Omega)}^{osc}<\infty\Big\},$$
    where
    $$||f||_{H^\alpha_p(\Omega)}^{osc}:=||f||_{L^p(\R)}+\big|\big|\big(\int_0^1t^{-2\alpha}(\osc^Mf(\cdot,t))^2\frac{dt}t\big)^{1/2}\,\big|\big|_{L^p(\R)}.$$
\end{theorem}
For a proof we refer to \cite[Thm. 3.5.1]{Tri92}.\\
We can now state and proof the extension theorem needed in the proof of Lemma \ref{Lemma1}.
\begin{theorem}[\cite{Tri92}, Theorem 4.5.2]
    Let $1<\alpha^{-1}<p<\infty,\, 0<\lambda_0<\lambda_1$ and $ a_0,a_1\in\R$ with 
    $$a_0+a_1=1,\quad a_0\lambda_0+a_1\lambda_1=-1.$$
    Then $\exti$ given by
    $$\exti f(x)=\begin{cases}
        f(x), &\text{ if }x\in \R_+,\\
        a_0f(-\lambda_0 x)+a_1f(-\lambda_1x),&\text{ else,}
    \end{cases}$$
    is a continuous extension operator from $H^\alpha_p(\R_+)$ to $H^\alpha_p(\R)$.
\end{theorem}
 
\begin{proof}
    In the following we will write $\osc^1 f$ for oscillations with regard to $\Omega=\R$ and $\Osc f$ for oscillations with regand to $\Omega=\R_+$, see Definition \ref{Osc}.\\
    Due to Lemma \ref{ptwsBess} $\exti$ makes sense for $f\in H^\alpha_p(\R_+)$.
    For $p\in\mathcal P_1$ one can easily verify, that $\exti p(x)=p(x)$ for all $x\in\R$. \\
    %In particular $\exti (f-p)=\exti (f)-p$ due to linearity for all $f\in H^\alpha_p(\R_+)$ and $p\in\mathcal P_1$.\\
    Define $a:=1+|a_0|+|a_1|$, $\lambda:=1+\lambda_0+\lambda_1$ and $\tilde \lambda:=\lambda/\min\{\lambda_0,\lambda_1\}$.
    Take now $0<t\leq 1$. Then we get for $x>t$
    $$\osc^1(\exti f)(x,t)=\Osc f(x,t)\leq \lambda\Osc f(x,\lambda t).$$ %
    For $x<-t$ we get due to $p=\exti(p|_{\R_+})$
    \begin{eqnarray*}
    \begin{aligned}
    \osc^1 (\exti f)(x,t)=&\inf_{p\in\mathcal P_1}\intbar_{B(x,t)}|\exti f(y)-p(y)|\,dy\\
   % =&\inf_{p\in\mathcal P_1}\intbar_{B(x,t)}|\exti f(y)-\exti p(y)|\,dy\\
    =&\inf_{p_0,p_1\in\mathcal P_1}\intbar_{B(x,t)}|\exti f(y)-a_0p_0(-\lambda_0y)-a_1p_1(-\lambda_1y)|\,dy\\
     =&\inf_{p_0,p_1\in\mathcal P_1}\intbar_{B(x,t)}\Big|\sum_{i=0}^1a_i(f(-\lambda_iy)-p_i(-\lambda_iy))\Big|\,dy.\\
    \end{aligned}
    \end{eqnarray*}
    These identities imply
     \begin{eqnarray*}
    \begin{aligned}
    &\osc^1 (\exti f)(x,t)\leq\sum_{i=0}^1\inf_{p_i\in\mathcal P_1}\intbar_{B(x,t)}|a_i(f(-\lambda_iy)-p_i(-\lambda_iy))|\,dy\\
     %\leq&a\sum_{i=0}^1\inf_{p_i\in\mathcal P_1}\oint_{B(x,t)}|f(-\lambda_iy)-p_i(-\lambda_iy)|\,dy\\
     \leq& a\sum_{i=0}^1\inf_{p_i\in\mathcal P_1}\frac1{2t}\int_{x-t}^{x+t}|f(-\lambda_iy)-p_i(-\lambda_iy)|\,dy
     =a\sum_{i=0}^1\inf_{p_i\in\mathcal P_1}\frac1{2t\lambda_i}\int_{\lambda_i(-x-t)}^{\lambda_i(-x+t)}|f(z)-p_i(z)|\,dz\\
     \leq& a\sum_{i=0}^1\inf_{p_i\in\mathcal P_1}\frac1{2t\lambda_i}\int_{(-\lambda_ix-\lambda t)_+}^{-\lambda_ix+\lambda t}|f(z)-p_i(z)|\,dz
     \leq a\tilde\lambda\sum_{i=0}^1\inf_{p_i\in\mathcal P_1}\frac1{2t\lambda}\int_{(-\lambda_ix-\lambda t)_+}^{-\lambda_ix+\lambda t}|f(z)-p_i(z)|\,dz\\
     \leq& a\tilde\lambda \sum_{i=0}^1\Osc f(-\lambda_i x,\lambda t).
    \end{aligned}
    \end{eqnarray*}
    For $-t\leq x\leq t$ we get with similar arguments
     \begin{eqnarray*}
    \begin{aligned}
    &\osc^1 (\exti f)(x,t)=\inf_{p\in\mathcal P_1}\intbar_{B(x,t)}|\exti f(y)-p(y)|\,dy\\
    %=&\inf_{p\in\mathcal P_1}\oint_{B(x,t)}|\exti f(y)-p(y)|\,dy\\
    =&\inf_{p\in\mathcal P_1}\frac{1}{2t}\Big(\int_{x-t}^0 |\sum_{i=0}^1a_i(f(-\lambda_iy)-p(-\lambda_iy))|\,dy +\int_0^{x+t}|f(y)-p(y)|\,dy\Big)\\
    \leq& \inf_{p\in\mathcal P_1}\frac{1}{2t}\Big(\frac{a\tilde \lambda}{\lambda}\sum_{i=0}^1 \int_{0}^{\lambda_i|x|+\lambda_i t} |f(z)-p(z)|\,dz +\int_0^{|x|+t}|f(y)-p(y)|\,dy\Big)\\
    %\leq& \inf_{p\in\mathcal P_1}\frac{\tilde \lambda|a|}{2t}\Big(\sum_{i=0}^1\int_{0}^{\lambda_ix+\lambda t} |f(y)-p(y)|\,dy +\int_0^{x+\lambda t}|f(y)-p(y)|\,dy\Big)\\
     \leq& \inf_{p\in\mathcal P_1}\frac{1}{2t}\Big(\frac{a\tilde \lambda}{\lambda}\sum_{i=0}^1\int_{0}^{\lambda |x|+\lambda t} |f(z)-p(z)|\,dz +\int_0^{\lambda |x|+\lambda t}|f(y)-p(y)|\,dy\Big)\\
      \leq& \inf_{p\in\mathcal P_1}\frac{1}{2t\lambda}(2a\tilde \lambda+\lambda)\int_{0}^{\lambda |x|+\lambda t} |f(y)-p(y)|\,dy
      \leq (2a\tilde \lambda+\lambda) \Osc f(\lambda |x|,\lambda t).\\
    \end{aligned}
    \end{eqnarray*}
   % and analogously for $-t\leq x<0$ the inequality
   % $$\osc^1 (\exti f)(x,t) \leq 3\tilde \lambda|a|\lambda \Osc^1 f(-\lambda x,\lambda t).$$
    We therefore get for $x\in \R$ and $0<t\leq 1$
    \begin{eqnarray*}
    \begin{aligned}
    \osc \exti f(x,t)\leq& (2a\tilde \lambda+\lambda) \Osc f(\lambda |x|,\lambda t)+a\tilde\lambda \sum_{i=0}^1\Osc f(\lambda_i |x|,\lambda t)+\lambda\Osc f(|x|,\lambda t)\\
    \leq& C\big(\Osc  f(\lambda |x|,\lambda t)+\sum_{i=0}^1\Osc  f(\lambda_i |x|,\lambda t)+\Osc f(|x|,\lambda t)\big) 
    \end{aligned}
    \end{eqnarray*}
    for some $C>0$. Now one gets for $f\in H^\alpha_p(\R_+)$
     \begin{eqnarray*}
    \begin{aligned}
    &||\exti f||^{osc}_{H^\alpha_p(\R)}\\
    =& ||\exti f||_{L^p(\R)}+\big|\big|\big(\int_0^1t^{-2\alpha}(\osc^1\exti f(\cdot,t))^2\frac{dt}t\big)^{1/2}\,\big|\big|_{L^p(\R)}\\
    \leq & C\Big(||f||_{L^p(\R_+)}+\\&\big|\big|\big(\int_0^1t^{-2\alpha}\big(\Osc  f(\lambda |x|,\lambda t)+\sum_{i=0}^1\Osc f(\lambda_i |x|,\lambda t)+\Osc f(|x|,\lambda t) \big)^2 \frac{dt}t\big)^{1/2}\,\big|\big|_{L^p(\R)}\Big)\\
    \leq & C\Big(||f||_{L^p(\R_+)}+\\&\big|\big|\big(\int_0^1t^{-2\alpha}\big((\Osc f(\lambda |x|,\lambda t))^2+\sum_{i=0}^1(\Osc f(\lambda_i |x|,\lambda t))^2+(\Osc f(|x|,\lambda t))^2\big)\frac{dt}t\big)^{1/2}\,\big|\big|_{L^p(\R)}\Big)\\
    \leq & C\Big(||f||_{L^p(\R_+)}+\\&\big|\big|\big(\int_0^1t^{-2\alpha}(\Osc f(\lambda |x|,\lambda t))^2\frac{dt}t\big)^{1/2}\,\big|\big|_{L^p(\R)} 
    +\big|\big|\big(\int_0^1t^{-2\alpha}(\Osc f(|x|,\lambda t))^2\frac{dt}t\big)^{1/2}\,\big|\big|_{L^p(\R)}\\
    +&\big|\big|\big(\int_0^1t^{-2\alpha}(\Osc f(\lambda_0 |x|,\lambda t))^2\frac{dt}t\big)^{1/2}\,\big|\big|_{L^p(\R)}+
    \big|\big|\big(\int_0^1t^{-2\alpha}(\Osc f(\lambda_1 |x|,\lambda t))^2\frac{dt}t\big)^{1/2}\,\big|\big|_{L^p(\R)}\Big)\\
    \leq& C\Big(||f||_{L^p(\R_+)}+ \big|\big|\big(\int_0^1t^{-2\alpha}(\Osc f(\lambda|x|,\lambda t))^2\frac{dt}t\big)^{1/2}\,\big|\big|_{L^p(\R)}\Big),
    \end{aligned}
    \end{eqnarray*}
    where we used linear substitutions in $x$ in the last inequality.\\
   Now observe for $\tilde f(\cdot):=f(\lambda\cdot)$, $ \tilde f\in H^\alpha_p(\R_+)$ and $x\in\R$, $t\in(0,1]$
    \begin{eqnarray*}
    \begin{aligned}
    \Osc f(\lambda |x|,\lambda t)\leq& \inf_{p\in P_1}\frac{1}{\lambda t}\int_{(\lambda |x|-\lambda t)_+}^{\lambda |x|+\lambda t} |f(y)-p(y)|\, dy\\
    =&\inf_{p\in P_1}\frac{1}{ t}\int_{(|x|- t)_+}^{|x|+ t} |f(\lambda z)-p(z)|\, dz\\
    \leq& 2\Osc \tilde f( |x|, t).
    \end{aligned}
    \end{eqnarray*}
    This leads to 
     \begin{eqnarray*}
    \begin{aligned}
    &||f||_{L^p(\R_+)}+ \big|\big|\big(\int_0^1t^{-2\alpha}(\Osc f(\lambda|x|,\lambda t))^2\frac{dt}t\big)^{1/2}\,\big|\big|_{L^p(\R)}\\
    \leq &C\Big(||\tilde f||_{L^p(\R_+)}+ \big|\big|\big(\int_0^1t^{-2\alpha}(\Osc \tilde f(|x|,t))^2\frac{dt}t\big)^{1/2}\,\big|\big|_{L^p(\R)}\Big)\\
     \leq &C\Big(||\tilde f||_{L^p(\R_+)}+ 2\big|\big|\big(\int_0^1t^{-2\alpha}(\Osc \tilde f(x,t))^2\frac{dt}t\big)^{1/2}\,\big|\big|_{L^p(\R_+)}\Big)\\
     \leq& C ||g||_{H^\alpha_p(\R)}^{osc}
    \end{aligned}
    \end{eqnarray*}
    for every extension $g\in H^\alpha_p(\R)$ of $\tilde f$. This implies that there is a constant $C>0$ independent of $f$ such that
    $$||\exti f||_{H^\alpha_p(\R)}\leq C ||\tilde f||_{H^\alpha_p(\R_+)}.$$
 By using the equivalent norm $||\cdot||_{L^p(\R)}+||\cdot||_{\mathring H^\alpha_p(\R)}$ on $H^\alpha_p(\R)$ one can quickly show that for $\eta\neq0$ there are constants $C_1,C_2>0$, such that
  $$||f(\cdot)||_{H^\alpha_p(\R)}\leq C_1||f(\eta \cdot)||_{H^\alpha_p(\R)}\leq C_2||f(\cdot)||_{H^\alpha_p(\R)}.$$
  The inequalities for the $L^p$-norm are proved easily. Furthermore
 %  For $f\in H^\alpha_p(\R)$ and $\eta\in \R\setminus\{0\}$ it is
\begin{eqnarray*}
\begin{aligned}
&||f(\eta\ \cdot)||_{\mathring H^{\alpha}_p(\R)}=||\F^{-1}\big(|\xi |^\alpha \F( f(\eta\ \cdot))(\xi)\big) ||_{L^p} \\
=& ||\F^{-1}\big(|\xi |^\alpha \frac{1}{|\eta|} (\F f)(\xi /\eta)\big) ||_{L^p}=
|\eta|^\alpha || \F^{-1}\big(\Big |\frac{\xi}{\eta}\Big |^\alpha \frac{1}{|\eta|} (\F f)(\xi /\eta)\big) ||_{L^p}\\
=&|\eta|^\alpha || (\F^{-1}(|\cdot |^\alpha (\F f)))(\eta\ \cdot) ||_{L^p} = |\eta|^{\alpha -1/p} || (\F^{-1}(|\cdot |^\alpha (\F f)))(\cdot) ||_{L^p}\\
=&|\eta|^{\alpha -1/p} ||f||_{\mathring H^{\alpha}_p(\R)}.
\end{aligned}
\end{eqnarray*}
  We immediately get analogous statements for the norms on $\R_+$. This implies
      \begin{eqnarray*}
    \begin{aligned}
    ||\exti f||_{H^\alpha_p(\R)}
   \leq C||\tilde f||_{H^\alpha_p(\R_+)}\leq C|| f||_{H^\alpha_p(\R_+)}.
    \end{aligned}
    \end{eqnarray*}
\end{proof}

\bibliographystyle{amsplain}
\bibliography{References}

\end{document}